\documentclass[11pt,reqno,a4paper]{amsart}
\usepackage{amssymb,amscd,amsmath,mathptmx}

\usepackage[matrix,arrow,curve,frame]{xy}    

\xymatrixcolsep{1.9pc}                          
\xymatrixrowsep{1.9pc}
\newdir{ >}{{}*!/-5pt/\dir{>}}                  

\makeatletter

\renewcommand{\subsection}{\@startsection{subsection}{1}{0pt}{-3.25ex plus -1ex minus-.2ex}{1.5ex plus.2ex}{\normalfont\it}}
\renewcommand{\section}{\@startsection{section}{1}{\parindent}{3.5ex plus 1ex minus .2ex}{2.3ex plus.2ex}{\sc}}

\renewcommand{\phi}{\varphi}
\renewcommand{\to}{\longrightarrow}
\renewcommand{\leq}{\leqslant}
\renewcommand{\geq}{\geqslant}
\renewcommand{\epsilon}{\varepsilon}

\renewcommand{\kappa}{\varkappa}

\DeclareMathOperator{\cyl}{cyl} \DeclareMathOperator{\Cyl}{Cyl}
\DeclareMathOperator{\hocolim}{hocolim}
 \DeclareMathOperator{\Map}{Map}
\DeclareMathOperator{\map}{map}

 \DeclareMathOperator{\Hom}{Hom}

\DeclareMathOperator{\id}{id} 
 
\DeclareMathOperator{\Alg}{Alg} \DeclareMathOperator{\colim}{colim}
\DeclareMathOperator{\Ho}{Ho}

 \DeclareMathOperator{\kr}{Ker}
 
 \DeclareMathOperator{\im}{Im}

\newcommand{\sn}{\par\smallskip\noindent}

\newcommand {\lp}{\varinjlim}

\newcommand{\lra}[1]{\bl{#1}\longrightarrow\relax}
\newcommand{\bl}[1]{\buildrel #1\over}
\newcommand{\cc}{\mathcal}
\newcommand{\bb}{\mathbb}
\newcommand{\ps}{\oplus}
\newcommand{\ff}{\mathfrak}

\newcommand{\op}{{\textrm{\rm op}}}

\newcommand{\wh}{\widehat}
\newcommand{\wt}{\widetilde}

\newcommand{\ifff}{if and only if }
\newcommand{\ring}{\cc Ring}

\newcommand{\aha}{\Alg_H}

\newtheorem{thm}{Theorem}[section]
\newtheorem{prop}[thm]{Proposition}
\newtheorem{cor}[thm]{Corollary}
\newtheorem{lem}[thm]{Lemma}
\newtheorem*{question}{Question}

\newtheorem*{rem}{Remark}
\newtheorem*{exs}{Examples}

\newtheorem*{defs}{Definition}
\newtheorem*{theo}{Theorem}

\makeatother

\begin{document}

\footskip30pt


\title{Homotopy theory of associative rings}
\author{Grigory Garkusha}
\address{School of Mathematics, University of Manchester, Oxford Road, M13~9PL Manchester, UK}
\urladdr{http://homotopy.nm.ru} \email{garkusha@imi.ras.ru}
\begin{abstract}
A kind of unstable homotopy theory on the category of associative
rings (without unit) is developed. There are the notions of
fibrations, homotopy (in the sense of Karoubi), path spaces, Puppe
sequences, etc. One introduces the notion of a quasi-isomorphism (or
weak equivalence) for rings and shows that - similar to spaces - the
derived category obtained by inverting the quasi-isomorphisms is
naturally left triangulated. Also, homology theories on rings are
studied. These must be homotopy invariant in the algebraic sense,
meet the Mayer-Vietoris property and plus some minor natural axioms.
To any functor $\cc X$ from rings to pointed simplicial sets a
homology theory is associated in a natural way. If $\cc X=GL$ and
fibrations are the $GL$-fibrations, one recovers
Karoubi-Villamayor's functors $KV_i, i>0$. If $\cc X$ is Quillen's
$K$-theory functor and fibrations are the surjective homomorphisms,
one recovers the (non-negative) homotopy $K$-theory in the sense of
Weibel. Technical tools we use are the homotopy information for the
category of simplicial functors on rings and the Bousfield
localization theory for model categories.  The machinery developed
in the paper also allows to give another definition for the
triangulated category $kk$ constructed by Corti\~nas and
Thom~\cite{CT}. The latter category is an algebraic analog for
triangulated structures on operator algebras used in Kasparov's
$KK$-theory.
\end{abstract}
\maketitle

\thispagestyle{empty} \pagestyle{plain}

\newdir{ >}{{}*!/-6pt/@{>}} 

\tableofcontents

\section{Introduction}

In the sixties mathematicians invented lower algebraic $K$-groups of
a ring and proved various exact sequences involving $K_0$ and $K_1$
(see Bass~\cite{Bass}). For instance, given a cartesian square of
rings
   \begin{equation}\label{222}
    \xymatrix{A\ar[r]\ar[d]&B\ar[d]^f\\
               C\ar[r]^g&D}
   \end{equation}
with $f$ or $g$ surjective, Milnor~\cite{Bass} proved a
Mayer-Vietoris sequence involving $K_0$ and $K_1$: the induced
sequence of abelian groups
   \begin{equation}\label{milnor}
    K_1(A)\xrightarrow{} K_1(C)\ps K_1(B)\xrightarrow{}
    K_1(D)\xrightarrow{\partial} K_0(A)\xrightarrow{} K_0(C)\ps
    K_0(B)\xrightarrow{} K_0(D)
   \end{equation}
is exact.

After Quillen~\cite{Q1} the higher algebraic $K$-groups of a ring
$R$ are defined by producing a space $K(R)$ and setting
$K_n(R)=\pi_nK(R)$. $K$ can be defined so that it actually gives a
functor $(Rings)\to (Spaces)$, and so the groups $K_n(R)$ start to
look like a homology theory on rings. However, there are negative
results which limit any search for extending the exact
sequence~\eqref{milnor} to the left involving higher $K$-groups. For
example, Swan~\cite{S} has shown that there is no satisfactory
$K$-theory, extending $K_0$ and $K_1$ and yielding Mayer-Vietoris
sequences, even if both $f$ and $g$ are surjective. Moreover, the
algebraic $K$-theory is not homotopy invariant in the algebraic
sense. These remarks show that $K$ is not a homology theory in the
usual sense.

Given any admissible category $\Re$ of rings with or without unit
(defined in Section~\ref{prel}) Gersten~\cite{G} considers group
valued functors $G$ on $\Re$ which preserve zero object, cartesian
squares, and kernels of surjective ring homomorphisms. He calls such
a functor a left exact MV-functor. It leads naturally to a homology
theory $\{k_i^G, i\geq 1\}$ of group valued functors on $\Re$. We require a
homology theory to be homotopy invariant in the algebraic sense,
to meet the Mayer-Vietoris property, and some other minor natural
properties given in Section~\ref{hom}. If
$G=GL$ one recovers the functors $KV_i$ of Karoubi and
Villamayor~\cite{KV}. The groups $KV_i(A)$ coincide with
$K_i(A)$ for any regular ring $A$.

After developing the general localization theory for model
categories in the 90-s (see the monograph by Hirschhorn~\cite{Hir})
we now have new devices for producing homology theories on rings.
More precisely, we fix an admissible category of rings $\Re$ and a
family of fibrations $\ff F$ on it like, for example, the
$GL$-fibrations or the surjective homomorphisms. Then any simplicial
functor on $\Re$ gives rise to a homology theory:

\begin{theo}
To any functor $\cc X$ from $\Re$ to pointed simplicial sets a
homology theory $\{k_i^{\cc X}, i\geq 0\}$ is associated. Such a
homology theory is defined by means of an explicitly constructed
functor $Ex_{I,J}(\cc X)$ from $\Re$ to pointed simplicial sets and,
by definition,
   $$k_i^{\cc X}(A):=\pi_i(Ex_{I,J}(\cc X)(A))$$
for any $A\in\Re$ and $i\geq 0$. Moreover, there is a natural
transformation $\theta_{\cc X}:\cc X\to Ex_{I,J}(\cc X)$,
functorial in $\cc X$.
\end{theo}

Roughly speaking, we turn any pointed simplicial functor into a
homology theory. If $\cc X=G$ one recovers the functors $k_i^G$ of
Gersten. In this way, the important simplicial functors $GL$ and $K$
give rise to the homology theories $\{KV_i\mid\ff
F=\textrm{$GL$-fibrations}\}$ and $\{KH_i\mid\ff
F=\textrm{surjective maps}\}$ respectively. Here $KH$ stands for the
(non-negative) homotopy $K$-theory in the sense of Weibel~\cite{W1}.

Next we present another part, developing a sort of unstable homotopy
theory on an admissible category of associative rings $\Re$. We are
based on the feeling that if rings are in a certain sense similar to
spaces then there should exist a homotopy theory where the
homomorphism $A\to A[x]$ is a homotopy equivalence, the Puppe
sequence, constructed by Gersten in~\cite{G1}, leads to various long
exact sequences, the loop ring $\Omega A=(x^2-x)A[x]$ is interpreted
as the loop space, etc.

For this we give definitions of quasi-isomorphisms for rings and
left derived categories $D^-(\Re,\ff F)$ associated to any family of
fibrations $\ff F$ on $\Re$. We show how to construct $D^-(\Re,\ff
F)$, mimicking the passage from spaces or chain complexes to the
homotopy category and the localization from this homotopy category
to the derived category.

In this way, the left derived category $D^-(\Re,\ff F)$ is obtained
from the admissible category of rings $\Re$ in two stages. First one
constructs a quotient $\cc H\Re$  of $\Re$ by equating homotopy
equivalent (in the sense of Karoubi) homomorphisms between rings.
Then one localizes $\cc H\Re$ by inverting quasi-isomorphisms via a
calculus of fractions. These steps are explained in
Section~\ref{der}. If $\ff F$ is saturated, which is always the case
in practice, then $D^-(\Re,\ff F)$ is naturally left triangulated.
The left triangulated structure as such is a tool for producing
homology theories on rings.

\begin{theo}
Let $\ff F$ be a saturated family of fibrations in $\Re$. One can
define the category of left triangles $\cc {L}tr(\Re,\ff F)$ in
$D^-(\Re,\ff F)$ having the usual set of morphisms from $\Omega
C\lra{f}A\lra{g}B\lra{h}C$ to $\Omega
C'\lra{f'}A'\lra{g'}B'\lra{h'}C'$. Then $\cc {L}tr(\Re,\ff F)$ is a
left triangulation of $D^-(\Re,\ff F)$, i.e. it is closed under
isomorphisms and enjoys the axioms which are versions of Vierdier's
axioms for triangulated categories. Stabilization of the loop
functor $\Omega$ produces a triangulated category $D(\Re,\ff F)$ out
of the left triangulated category $D^-(\Re,\ff F)$.
\end{theo}

Motivated by ideas and work of J. Cuntz on bivariant $K$-theory of
locally convex algebras (see~\cite{Cu,Cu1}), Corti\~nas and Thom~\cite{CT} construct a
bivariant homology theory $kk_*(A,B)$ on the category $\aha$ of
algebras over a unital ground ring $H$. It is Morita invariant,
homotopy invariant, excisive $K$-theory of algebras, which is
universal in the sense that it maps uniquely to any other such
theory. This bivariant $K$-theory is defined in a triangulated
category $kk$ whose objects are the $H$-algebras without unit and
$kk_n(A,B)=kk(A,\Omega^{n}B)$, $n\in\bb Z$. We make use of our
machinery to study various triangulated structures on admissible
categories of rings which are not necessarily small. As an
application, we give another, but equivalent, description of the
triangulated category $kk$.

\begin{theo}
Let $\Re$ be an arbitrary admissible category of rings and let $\ff
W$ be any subcategory of homomorphisms containing $A\to A[x]$ such
that the triple $(\Re,\ff W, \ff F=\{\textrm{surjective maps}\})$ is
a Brown category. There is a triangulated category $D(\Re,\ff W)$
whose objects and morphisms are defined similar to those of
$D(\Re,\ff F)$. If $\Re=\aha$ and $\ff W_{CT}$ is the class of weak
equivalences generated by Morita invariant, homotopy invariant,
excisive homology theories, then there is a natural triangulated
equivalence of the triangulated categories $D(\aha,\ff W_{CT})$ and
$kk$.
\end{theo}

The main tools of the paper are coming from modern homotopical
algebra (as exposed for instance in the work of Hovey~\cite{H},
Hirschhorn~\cite{Hir}, Dugger~\cite{D}, Goerss and
Jardine~\cite{GJ}). To develop homotopy theory of rings we consider
the model category $U\Re$ of simplicial functors on $\Re$, i.e.
simplicial presheaves on $\Re^{\op}$ instead of simplicial
presheaves on $\Re$. The model structure is given by injective maps
(cofibrations) and objectwise weak equivalences of simplicial sets
(Quillen equivalences). There is a contravariant embedding $r$ of
$\Re$ into $U\Re$ as representable functors. We need to localize
this model structure to take into account the pullback
squares~\eqref{222} with $f$ a fibration in $\ff F$ and the fact
that $rA[x]\to rA$ should be a Quillen equivalence. Let us remark
that we require a homology theory to take such distinguished squares
to the Mayer-Vietoris sequence. To do so, we define a set $\cc S$ to
consist of the maps $rA[x]\to rA$ for any ring $A$ and maps
$rB\bigsqcup_{rD}rC\to rA$ for every pullback square~\eqref{222} in
$\Re$ with $f$ a fibration. Then one localizes $U\Re$ at $\cc S$.
This procedure is a reminiscence of an unstable motivic model
category. The latter model structure is obtained from simplicial
presheaves $\cc E$ on smooth schemes by localizing $\cc E$ at the
set $\cc S$ of the maps $X\times\bb A^1\to X$ for any smooth scheme
$X$ and maps $P\to D$ for every pullback square~\eqref{222} of
smooth schemes with $f$ etale, $g$ an open embedding, and
$f^{-1}(D-C)\to D-C$ an isomorphism. There is then some work
involving properties of the Nisnevich topology to show that this
model category is equivalent to the Morel-Voevodsky motivic model
category of~\cite{MV}.

\subsubsection*{Organization of the paper} After fixing some notation and
terminology in Section~\ref{prel}, we study the notion of
$I$-homotopy for simplicial functors on an admissible category of
rings $\Re$. It has a lot of common properties with $\bb
A^1$-homotopy for simplicial (pre-)sheaves on schemes. We show there
how to convert a simplicial functor into a homotopy invariant one.
All this material is the content of Section~\ref{sin}. Then comes
Section~\ref{hom} in which homology theories on rings are
investigated. We also construct there the simplicial functor
$Ex_{I,J}(\cc X)$. Derived categories on rings and their left
triangulated structure are studied in Section~\ref{der}. In
Section~\ref{stab} the stabilization procedure is described as well
as the triangulated categories $D(\Re,\ff F)$. In Section~\ref{kk}
we apply the machinery developed in the preceding sections to study
various triangulated structures on admissible categories of rings
which are not necessarily small. We also give an equivalent
definition of $kk$ there. The necessary facts about Bousfield
localization in model categories are given in Addendum.

\subsubsection*{Acknowledgement} This paper was written during the visits of the author in `05 to
the Euler IMI in St.~Petersburg and IHES in Paris and completed
during the visit in `06 to the University of Manchester (supported
by the MODNET Research Training Network in Model Theory). He would
like to thank the Institutes and the University for the kind
hospitality.

\section{Preliminaries}\label{prel}

We shall work in the category $\ring$ of associative rings (with
or without unit) and ring homomorphisms. Following
Gersten~\cite{G} a category of rings $\Re$ is {\it admissible\/}
if it is a full subcategory of $\ring$ and

\begin{enumerate}

\item $R$ in $\Re$, $I$ a (two-sided) ideal of $R$ then $I$ and
$R/I$ are in $\Re$;

\item if $R$ is in $\Re$, then so is $R[x]$, the polynomial ring in
one variable;

\item given a cartesian square
   $$\xymatrix{D\ar[r]^\rho\ar[d]_\sigma &A\ar[d]^f\\
               B\ar[r]^g &C}$$
in $\ring$ with $A,B,C$ in $\Re$, then $D$ is in $\Re$.

\end{enumerate}

One may abbreviate 1, 2, and 3 by saying that $\Re$ is closed under
operations of taking ideals, homomorphic images, polynomial
extensions in a finite number of variables, and fibre products. If
otherwise stated we shall always work in a fixed (skeletally) small
admissible category $\Re$.

\begin{rem}{\rm
Given a ring homomorphism $f:R\to R'$ in $\ring$ between two rings
with unit, $f(1)$ need not be equal to $1$. We only assume that
$f(r_1r_2)=f(r_1)f(r_2)$ and $f(r_1+r_2)=f(r_1)+f(r_2)$ for any two
elements $r_1,r_2\in R$. It follows that the trivial ring 0 is a
zero object in $\ring$. }\end{rem}

If $R$ is a ring then the polynomial ring $R[x]$ admits two
homomorphisms onto $R$
   $$\xymatrix{R[x]\ar@<2.5pt>[r]^{\partial_x^0}\ar@<-2.5pt>[r]_{\partial_x^1}&R}$$
where
   $$\partial_x^i|_R=1_R,\ \ \ \partial_x^i(x)=i,\ \ \ i=0,1.$$
Of course, $\partial_x^1(x)=1$ has to be understood in the sense
that $\Sigma r_nx^n\mapsto\Sigma r_n$.

\begin{defs}{\rm
Two ring homomorphisms $f_0,f_1:S\to R$ are {\it elementary
homo\-topic}, written $f_0\sim f_1$, if there exists a ring
homomorphism
   $$f:S\to R[x]$$
such that $\partial^0_xf=f_0$ and $\partial^1_xf=f_1$. A map
$f:S\to R$ is called an {\it elementary homotopy equivalence\/} if
there is a map $g:R\to S$ such that $fg$ and $gf$ are elementary
homotopic to $\id_R$ and $\id_S$ respectively.

}\end{defs}

For example, let $A$ be a $\bb N$-graded ring, then the inclusion
$A_0\to A$ is an elementary homotopy equivalence. The homotopy
inverse is given by the projection $A\to A_0$. Indeed, the map $A\to
A[x]$ sending a homogeneous element $a_n\in A_n$ to $a_nt^n$ is a
homotopy between the composite $A\to A_0\to A$ and the identity
$\id_A$.

The relation ``elementary homotopic'' is reflexive and
symmetric~\cite[p.~62]{G}. One may take the transitive closure of
this relation to get an equivalence relation (denoted by the
symbol ``$\simeq$''). The set of equivalence classes of morphisms
$R\to S$ is written $[R,S]$.

\begin{lem}[Gersten \cite{G1}]
Given morphisms in $\ring$
   $$\xymatrix{R\ar[r]^f &S\ar@/^/[r]^g \ar@/_/[r]_{g'} &T\ar[r]^h &U}$$
such that $g\simeq g'$, then $gf\simeq g'f$ and $hg\simeq hg'$.
\end{lem}

Thus homotopy behaves well with respect to composition and we have
category $Hotring$, the {\it homotopy category of rings}, whose
objects are rings and such that $Hotring(R,S)=[R,S]$. The homotopy
category of an admissible category of rings $\Re$ will be denoted
by $\cc H(\Re)$.

The diagram in $\ring$
   $$A\bl{f}\to B\bl{g}\to C$$
is a short exact sequence if $f$ is injective ($\equiv \kr f=0$),
$g$ is surjective, and the image of $f$ is equal to the kernel of
$g$. Thus $f$ is a normal monomorphism in $\Re$ and $f=\ker g$.

\begin{defs}{\rm
A ring $R$ is {\it contractible\/} if $0\sim 1$; that is, if there
is a ring homomorphism $f:R\to R[x]$ such that $\partial^0_xf=0$
and $\partial^1_xf=1_R$.

}\end{defs}

Following Karoubi and Villamayor~\cite{KV} we define $ER$, the
{\it path ring\/} on $R$, as the kernel of $\partial_x^0:R[x]\to
R$, so $ER\to R[x]\bl{\partial_x^0}\to R$ is a short exact
sequence in $\ring$. Also $\partial_x^1:R[x]\to R$ induces a
surjection
   $$\partial_x^1:ER\to R$$
and we define the {\it loop ring\/} $\Omega R$ of $R$ to be its
kernel, so we have a short exact sequence in $\ring$
   $$\Omega R\to ER\bl{\partial_x^1}\to R.$$
Clearly, $\Omega R$ is the intersection of the kernels of
$\partial_x^0$ and $\partial_x^1$. By~\cite[3.3]{G} $ER$ is
contractible for any ring $R$.

\section{The functor $Sing_*$}\label{sin}

In this section we introduce and study the important notion of
$I$-homotopy for simplicial functors on an admissible category of
rings $\Re$. It is similar to $\bb A^1$-homotopy in the sense of
Morel and Voevod\-sky~\cite{MV}.

\subsection{Homotopization}

Recall that a simplicial set map $f:X\to Y$ is a weak equivalence
if all maps
\begin{enumerate}
\item $\pi_0X\to\pi_0 Y$, and
\item $\pi_i(X,x)\to\pi_i(Y,fx)$, $x\in X_0$, $i\geq 1$
\end{enumerate}
are bijections. Here $\pi_i(X,x)=\pi_i(|X|,x)$, in general, but
   $$\pi_i(X,x)=[(S^i,*),(X,x)]=\pi((S^i,*),(X,x))$$
if $X$ is fibrant (recall that
$S^i=\varDelta^i/\partial\varDelta^i$ is the simplicial
$i$-sphere).

Following Gersten, we say that a functor $F$ from rings to sets is
{\it homotopy invariant\/} if $F(R)\cong F(R[t])$ for every $R$.
Similarly, a functor $F$ from rings to simplicial sets is {\it
homotopy invariant\/} if for every ring $R$ the natural map $R\to
R[t]$ induces a weak equivalence of simplicial sets $F(R)\simeq
F(R[t])$. Note that each homotopy group $\pi_n(F(R))$ also forms a
homotopy invariant functor.

We shall introduce the simplicial ring $R[\varDelta]$, and use it
to define the homo\-to\-pi\-za\-tion functor $Sing_*$.

For each ring $R$ one defines a simplicial ring $R[\varDelta]$,
   $$R[\varDelta]_n:=R[\varDelta^n]=R[t_0,\ldots,t_n]/(\sum t_i-1)R\ \ \ (\cong
     R[t_1,\ldots,t_n]).$$
The face and degeneracy operators $\partial_i:R[\varDelta^n]\to
R[\varDelta^{n-1}]$ and $s_i:R[\varDelta^n]\to R[\varDelta^{n+1}]$
are given by
   $$\partial_i(t_j)\ (\textrm{resp. $s_i(t_j)$})=
     \left\{
      \begin{array}{rcl}
       t_j\ (\textrm{resp. $t_j$}),\ j&<&i\\
       0\ (\textrm{resp. $t_j+t_{j+1}$}),\ j&=&i\\
       t_{j-1}\ (\textrm{resp. $t_{j+1}$}),\ i&<&j
      \end{array}
      \right.$$
Note that the face maps $\partial_{0;1}:R[\varDelta^1]\to
R[\varDelta^{0}]$ are isomorphic to $\partial^{0;1}_t:R[t]\to R$
in the sense that the diagram
   $$\xymatrix{R[t]\ar[r]^{\partial^{\epsilon}_t}\ar[d]_{t\mapsto t_0}&R\ar[d]\\
               R[\varDelta^1]\ar[r]^{\partial_{\epsilon}}&R[\varDelta^{0}]}$$
is commutative and the vertical maps are isomorpisms.

\begin{lem}\label{nnn}
The inclusion of simplicial rings $R[\varDelta]\subset
R[x][\varDelta]$ is a homotopy equivalence, split by evaluation at
$x=0$.
\end{lem}

\begin{proof}
A simplicial homotopy from $R[x][\varDelta]$ to $R[x][\varDelta]$
is a simplicial map
   $$h:R[x][\varDelta]\times\varDelta^1\to R[x][\varDelta].$$
Recall that a $n$-simplex $v$ of $\varDelta^1$ is nothing more
than to give an integer $i$ with $-1\leq i\leq n$, and send the
integers $\{0,1,\ldots,i\}$ to 0, while the integers
$\{i+1,i+2,\ldots,n\}$ map to 1. So any homotopy is given by maps
   $$h^{(n)}_v:R[x][\varDelta^n]\to R[x][\varDelta^n],\ \ \ v\in\varDelta^1,$$
which must be compatible with the face and degeneracy operators.

Given $v=v(i)\in\varDelta^1$ let $h^{(n)}_v(f)=f$ if $f\in
R[\varDelta^n]$ and
   $$x\longmapsto
     \left\{
      \begin{array}{rcl}
       x(t_0+\cdots+t_i),\ i&\geq&0\\
       0\hskip1cm,\ i&=&-1
      \end{array}
      \right.$$
It is directly verified that the maps $h^{(n)}_v$ are compatible
with the face and degeneracy operators. These maps define a
simplicial homotopy between the identity map of $R[x][\varDelta]$
and the composite
   $$R[x][\varDelta]\xrightarrow{x=0}R[\varDelta]\subset R[x][\varDelta].$$
This implies the claim.
\end{proof}

\begin{defs}[Homotopization]{\rm
Let $F$ be a functor from rings to simplicial sets. Its {\it
homotopization\/} $Sing_*(F)$ is defined at each ring $R$ as the
diagonal of the bisimplicial set $F(R[\varDelta])$. Thus
$Sing_*(F)$ is also a functor from rings to simplicial sets. If we
consider $R$ as a constant simplicial ring, the natural map $R\to
R[\varDelta]$ yields a natural transformation $F\to Sing_*(F)$.

(Strict Homotopization). Let $F$ be a functor from rings to sets.
Its {\it strict homotopization\/} $[F]$ is defined as the
coequalizer of the evaluations at $t=0,1:F(R[t])\rightrightarrows
F(R)$. The coequaliser can be constructed as follows. Given
$x,y\in F(R)$, write $x\sim y$ if there is a $z\in F(R[t])$ such
that $(t=0)(z)=x$ and $(t=1)(z)=y$. Then this relation is
reflexive and symmetric (use the automorphism
$R[t]\xrightarrow{t\mapsto 1-t} R[t]$). Its transitive closure
determines an equivalence relation and then $[F](R)$ is the
quotient of $F(R)$ with respect to this equivalence relation.

In fact, $[F]$ is a homotopy invariant functor and there is a
universal transformation $F(R)\to[F](R)$. Moreover, if $F$ takes
values in groups then so does $[F]$ (see Weibel~\cite{W}).

}\end{defs}

Given a functor $F$ from rings to simplicial sets, by $F[t]$ denote
the functor which is defined as $F(R[t])$ at each ring $R$. The
natural inclusion $R\to R[t]$ yields a natural transformation $F\to
F[t]$.

\begin{prop}\label{ggg}
Let $F$ be a functor from rings to simplicial sets. Then:
\begin{enumerate}

\item $Sing_*(F)$ is a homotopy invariant functor;

\item if $F$ is homotopy invariant then $F(R)\to Sing_*(F)(R)$ is a weak equivalence for
all $R$ and $Sing_*(F)\to Sing_*(F)[t]$ is an objectwise homotopy
equivalence, functorial in $R$.

\item $\pi_0(Sing_*(F))$ is a strict homotopization $[F_0]$ of the functor
$F_0(R)=\pi_0(F(R))$.
\end{enumerate}
\end{prop}

\begin{proof}
Let us show that the inclusion of simplicial rings
$R[\varDelta]\subset R[x][\varDelta]$ induces a weak equivalence
$Sing_*(F)(R)\to Sing_*(F)(R[x])$. Actually we shall prove even
more: the latter map turns out to be a homotopy equivalence of
simplicial sets (also showing that $Sing_*(F)\to Sing_*(F)[t]$ is
a homotopy equivalence).

Let
   $$h^{(n)}_v:R[x][\varDelta^n]\to R[x][\varDelta^n],\ \ \ v\in\varDelta^1,$$
be the maps constructed in the proof of Lemma~\ref{nnn}. We claim
that the maps
   $$H^{(n)}_v=F_n(h^{(n)}_v):F_n(R[x][\varDelta^n])\to F_n(R[x][\varDelta^n]),\ \ \ v\in\varDelta^1,$$
define a simplicial homotopy between the identity map of
$Sing_*(F)(R[x])$ and the composite
   $$Sing_*(F)(R[x])\xrightarrow{x=0}Sing_*(F)(R)\to Sing_*(F)(R[x]).$$
For this, we must verify that the maps $H^{(n)}_v$ are compatible
with the structure maps $w:[m]\to[n]$ in $\varDelta$.

We already know that $w^*\circ h^{(n)}_v=h^{(m)}_{vw}\circ w^*$.
One has a commutative diagram
   $$\xymatrix{F_n(R[x][\varDelta^n])\ar[r]^{F_n(w^*)}\ar[d]_{w^*_F}\ar[dr]^(.6){w^*_{Sing_*(F)}}
               &F_n(R[x][\varDelta^m])\ar[rr]^{F_n(h^{(m)}_{vw})}\ar[d]^{w^*_F}&& F_n(R[x][\varDelta^m])\ar[d]^{w^*_F}\\
               F_m(R[x][\varDelta^n])\ar[r]_{F_m(w^*)}
               &F_m(R[x][\varDelta^m])\ar[rr]_{F_m(h^{(m)}_{vw})=H^{(m)}_{vw}} &&F_m(R[x][\varDelta^m]).}$$
Then,
   \begin{gather*}
    w^*_{Sing_*(F)}\circ H^{(n)}_v=w^*_F\circ F_n(w^*)\circ
    F_n(h^{(n)}_v)=w^*_F\circ F_n(h^{(m)}_{vw})\circ F_n(w^*)=\\
    H^{(m)}_{vw}\circ w^*_F\circ F_n(w^*)=H^{(m)}_{vw}\circ w^*_{Sing_*(F)}.
   \end{gather*}
We have checked that the maps $H^{(n)}_v$ are compatible with the
structure maps in $\varDelta$, as claimed. These give the
necessary simplicial homotopy. 

Part (3) follows from the fact that, for any simplicial space
$X.$, the group $\pi_0(|X.|)$ is the coequaliser of
$\partial_0,\partial_1:\pi_0(X_1)\rightrightarrows\pi_0(X_0)$. In
this case $\pi_0(X_0)=\pi_0(F(R))$ and
$\pi_0(X_1)=\pi_0(F(R[t]))$.
\end{proof}

Let $\Re$ be an admissible category of rings. In order to
construct homology theories on $\Re$, we shall use the model
category $U\Re$ of covariant functors from $\Re$ to simplicial
sets (and not contravariant functors as usual). Note that this
usage deviates from the usual notation and practice, e.g. as in
Dugger~\cite{D}. We consider the Heller model structure on $U\Re$
instead of the most commonly used Bousfield-Kan model structure.
It is a proper, simplicial, cellular model category with weak
equivalences and cofibrations being defined objectwise, and
fibrations being those maps having the right lifting property with
respect to trivial cofibrations (see Dugger~\cite{D}). We consider
the fully faithful contravariant functor
   $$r:\Re\to U\Re,\ \ \ A\longmapsto\Hom_{\Re}(A,-),$$
where $rA(B)=\Hom_{\Re}(A,B)$ is to be thought of as the constant
simplicial set for any $B\in\Re$.

The model structure on $U\Re$ enjoys the following properties (see
Dugger~\cite[p.~21]{D}):

\begin{itemize}
\item[$\diamond$] every object is cofibrant;

\item[$\diamond$] being fibrant implies being objectwise
fibrant, but is stronger (there are additional diagramatic
conditions involving maps being fibrations, etc.);

\item[$\diamond$] any object which is constant in the simplicial
direction is fibrant.

\end{itemize}
If $F\in U\Re$ then $U\Re(rA\times\varDelta^n,F)=F_n(A)$
(isomorphism of sets). Hence, if we look at simplicial mapping
spaces we find
   $$\Map(rA,F)=F(A)$$
(isomorphism of simplicial sets). This is a kind of ``simplicial
Yoneda Lemma''.

\begin{defs}{\rm
Let $f,g:\cc X\to\cc Y$ be two maps of simplicial presheaves in
$U\Re$. An {\it elementary $I$-homotopy\/} from $f$ to $g$ is a
map $H:\cc X\to\cc Y[t]$ such that $\partial^0\circ H=f$ and
$\partial^1\circ H=g$, where $\partial^{0;1}:\cc Y[t]\to\cc Y$ are
the maps induced by $\partial^{0;1}_t:A[t]\to A$, $A\in\Re$. Two
morphisms are said to be {\it $I$-homotopic\/} if they can be
connected by a sequence of elementary $I$-homotopies. A map $f:\cc
X\to\cc Y$ is called an {\it $I$-homotopy equivalence\/} if there
is a map $g:\cc Y\to\cc X$ such that $fg$ and $gf$ are
$I$-homotopic to $\id_{\cc Y}$ and $\id_{\cc X}$ respectively.

}\end{defs}

Let $A,B$ be two rings in $\Re$ and let $H:rB\to (rA)[t]$ be an
elementary homotopy of representable functors. It follows that the
map $\wh H:=H_B(\id_B):A\to B[t]$ yields an elementary homotopy
between $A$ and $B$. Moreover, for any ring $R\in\Re$ and any ring
homomorphism $\alpha:B\to R$
   \begin{equation}\label{aaa}
     H_R\circ\alpha_*(\id_B)=H_R(\alpha)=\alpha[t]\circ\wh H,
   \end{equation}
where $\alpha_*=\Re(B,\alpha):\Re(B,B)\to\Re(B,R)$ and
$\alpha[t]:B[t]\to R[t]$, $\sum b_it^i\longmapsto\sum\alpha
(b_i)t^i$.

Conversely, suppose $\wh H:A\to B[t]$ is an elementary homotopy in
$\Re$, then the collection of maps
   $$\{H_R(\alpha):=\alpha[t]\circ\wh H\mid R\in\Re,\alpha\in\Re(B,R)\}$$
gives rise to an elementary homotopy $H:rB\to (rA)[t]$.

\begin{cor}\label{qqq}
Two maps $f,g:A\to B$ are elementary homotopic in $\Re$ \ifff the
induced maps $f^*,g^*:rB\to(rA)[t]$ are elementary $I$-homotopic
in $U\Re$. Furthermore, there is a bijection between elementary
homotopies in $\Re$ and elementary $I$-homotopies in $U\Re$. This
bijection is given by~\eqref{aaa}.
\end{cor}

\subsection{The model category $U\Re_I$}

Let $\Re$ be an admissible category of rings and let
$I=\{i=i_A:r(A[t])\to r(A)\mid A\in\Re\}$, where each $i_A$ is
induced by the natural homomorphism $i:A\to A[t]$. We shall refer
to the $I$-local equivalences as $I$-weak equivalences. The
resulting model category $U\Re/I$ will be denoted by $U\Re_I$ and
its homotopy category is denoted by $\Ho_I(\Re)$. Notice that any
homotopy invariant functor $F:\Re\to Sets$ is an $I$-local object
in $U\Re$ (hence fibrant in $U\Re_I$).

The following lemma is straightforward.

\begin{lem}\label{www}
A fibrant object $\cc X\in U\Re$ is $I$-local \ifff the map $\cc
X\to\cc X[t]$ is a weak equivalence in $U\Re$.
\end{lem}

\begin{lem}
If two maps $f,g:\cc X\to\cc Y$ in $U\Re$ are elementary
$I$-homoto\-pic, then they coincide in the $I$-homotopy category
$\Ho_I(\Re)$.
\end{lem}

\begin{proof}
By assumption there is a map $H:\cc X\to\cc Y[t]$ such that
$\partial^0H=f$ and $\partial^1H=g$.

Let $\alpha:\cc Y\to\wh{\cc Y}$ be a fibrant replacement of $\cc
Y$ in $U\Re_I$. It follows that $\wh{\cc Y}$ is an $I$-local
object in $U\Re$. By Lemma~\ref{www} the map $i:\wh{\cc
Y}\to\wh{\cc Y}[t]$ is a weak equivalence.

One has a commutative diagram
   $$\xymatrix{\wh{\cc Y}\ar[r]^(.40){diag}\ar[d]_i&\wh{\cc Y}\times\wh{\cc Y}\\
               \wh{\cc Y}[t].\ar[ur]_{(\partial^0,\partial^1)}}$$
We see that $\wh{\cc Y}[t]$ is a path object of $\wh{\cc Y}$ in
$U\Re_I$.

Consider the following diagram:
   $$\xymatrix{&&&{\cc Y}[t]\ar@<2.5pt>[d]_{\partial^0\ }\ar@<-2pt>[d]^{\ \,\partial^1}\ar[rr]^{\alpha[t]}
               && \wh{\cc Y}[t]\ar@<2.5pt>[d]_{\partial^0\ }\ar@<-2pt>[d]^{\ \,\partial^1}\\
               &\cc X\ar@/^1.5pc/[urr]^H \ar@<2.5pt>[rr]^{f}\ar@<-2.5pt>[rr]_{g}
               && {\cc Y}\ar[rr]_{\alpha} && \wh{\cc Y}.}$$
Here
$\partial^\epsilon\circ\alpha[t]=\alpha\circ\partial^\epsilon$.
Hence $\alpha f\bl r\sim\alpha g$. It follows from~\cite[9.5.24;
9.5.15]{Hir} that $\alpha f$ and $\alpha g$ represent the same map
in the homotopy category. Since $\alpha$ is an isomorphism in
$\Ho_I(\Re)$, we deduce that $f=g$ in $\Ho_I(\Re)$.
\end{proof}

\begin{lem}\label{mmm}
Any $I$-homotopy equivalence is an $I$-weak equivalence.
\end{lem}

\begin{proof}
Let $f:\cc X\to\cc Y$ be an $I$-homotopy equivalence and $g$ be an
$I$-homotopy inverse to $f$. We have to show that the compositions
$fg$ and $gf$ are equal to the corresponding identity morphisms in
the $I$-homotopy category $\Ho_I(\Re)$. By definition, these maps
are $I$-homotopic to the identity and it remains to show that two
elementary $I$-homotopic morphisms coincide in the $I$-homotopy
category. But this follows immediately from the preceding lemma.
\end{proof}

\begin{lem}\label{eee}
For any $\cc X$ the canonical morphism $\cc X\to\cc X[t]$ is an
$I$-homotopy equivalence, and thus an $I$-weak equivalence.
\end{lem}

\begin{proof}
For any ring $R$ the natural homomorphism $i:R\to R[t]$ is an
elementary homotopy equivalence, split by evaluation at $t=0$.
Indeed, the homomorphism $R[t]\to R[t,y]$ sending $t$ to $ty$
defines an elementary homotopy between the identity homomorphism
and the composite
   $$R[t]\xrightarrow{t=0}R\subset R[t].$$
Applying $\cc X$ to the elementary homotopy equivalence $i:R\to
R[t]$, one gets an $I$-homoto\-py from $\cc X(i\circ(t=0))$ and
$\id_{\cc X[t]}$. Since $\cc X((t=0)\circ i)=\id_{\cc X}$, the
lemma is proven.
\end{proof}

\begin{cor}\label{fff}
For any $\cc X$ the canonical morphism $\cc X\to Sing_*(\cc X)$ is
an $I$-trivial cofibration.
\end{cor}

\begin{proof}
Since $R$ is a retract of $R[\varDelta^n]$ for any ring $R$ the
map of the assertion is plainly a cofibration. It remains to check
that it is an $I$-weak equivalence.

Given a functor $F:\Re\to Sets$, the canonical morphism $F\to
F[t_1,\ldots,t_n]$ is an $I$-weak equivalence by Lemma~\ref{eee}.
Since for any ring $R$ and any $n\geq 0$ the ring $R[\varDelta^n]$
is isomorphic to $R[t_1,\ldots,t_n]$, functorially in $R$, we see
that the canonical morphism $F\to F[\varDelta^n]$ is an $I$-weak
equivalence, where $F[\varDelta^n](R):=F(R[\varDelta^n])$.

The canonical morphism $\cc X\to Sing_*(\cc X)$ coincides
objectwise with the canonical morphisms $\cc X_n\to\cc
X_n[\varDelta^n]$. It follows from~\cite[18.5.3]{Hir} that the map
   $$\hocolim_{\varDelta^{\op}}\cc X_n\to\hocolim_{\varDelta^{\op}}\cc X_n[\varDelta^n]$$
is an $I$-weak equivalence. By~\cite[18.7.5]{Hir} the canonical
map $\hocolim_{\varDelta^{\op}}\cc X_n\to\cc X$ (respectively
$\hocolim_{\varDelta^{\op}}\cc X_n[\varDelta^n]\to Sing_*(\cc X)$)
is a weak equivalence in $U\Re$, whence the assertion follows.
\end{proof}

Let $\vartheta_{\cc X}:\cc X\to R({\cc X})$ denote a fibrant
replacement functor in $U\Re$. That is $R({\cc X})$ is fibrant and
the map $\vartheta_{\cc X}$ is a trivial cofibration in $U\Re$.
Given a model category $\cc C$, we write $\cc C_\bullet$ to denote
the model category under the terminal object~\cite[p. 4]{H}. If
$\cc C=U\Re$ we shall refer to the objects of $U\Re_\bullet$ as
pointed simplicial functors.

\begin{thm}\label{ccc}
The map $\cc X\longmapsto R({Sing_*(\cc X)})$ yields a fibrant
replacement functor in $U\Re_I$. That is the object $R({Sing_*(\cc
X)})$ is $I$-local and the composition
   $$\cc X\to {Sing_*(\cc X)}\to R({Sing_*(\cc X)})$$
is an $I$-trivial cofibration. Furthermore, the natural map
   $$\pi_0(Sing_*(\cc X)(A))=[\cc X_0](A)\to\Hom_{\Ho_I(\Re)}(rA,\cc X)$$
is a bijection for any $A\in\Re$. Moreover, if $\cc X$ is pointed,
then for any integer $n\geq 0$ and any $A\in\Re$ the obvious map
   $$\pi_n(Sing_*(\cc X)(A))\to\Hom_{\Ho_{I,\bullet}(\Re)}((rA_+)\wedge S^n,\cc X)$$
is a bijection, where $rA_+=rA\sqcup pt$.
\end{thm}

\begin{proof}
The fact that $R({Sing_*(\cc X)})$ is an $I$-local object is a
consequence of Proposition~\ref{ggg}. The map $\cc X\longmapsto
R({Sing_*(\cc X)})$ yields a fibrant replacement functor by
Corollary~\ref{fff}.

The rest of the proof follows from the fact that for any $\cc X\in
U\Re$ the function space of maps $\Map(rA,R(\cc X))$ may be
identified with $R(\cc X)(A)$, which is weakly equivalent to $\cc
X(R)$ because $\cc X\to R(\cc X)$ is an objectwise weak
equivalence.
\end{proof}

\begin{cor}\label{ddd}
Let $\Re$ be an admissible category of rings and let $\cc H\Re$ be
its homotopy category. Then the functor
   $$r:\cc H\Re\to\Ho_I(\Re),\ \ \ [A,B]\longmapsto\Ho_I(\Re)(rB,rA)$$
is a fully faithful contravariant embedding.
\end{cor}

\begin{proof}
This is a consequence of Proposition~\ref{ggg},
Corollary~\ref{qqq}, and Theorem~\ref{ccc}.
\end{proof}

Call a ring homomorphism $s:A\to B$ an {\it $I$-weak
equivalence\/} if its image in $U\Re$ is an $I$-weak equivalence.

\begin{cor}\label{dd}
Let $B$ be a ring in $\Re$ and consider a ring $B^I$ together with
homomorphisms
   $$B\lra{s}B^I\xrightarrow{(d_0,d_1)}B\times B,$$
where $s$ is an $I$-weak equivalence and the composite is the
diagonal. Then for any homomorphism $H:A\to B^I$ the homomorphisms
$d_0\circ H$ and $d_1\circ H$ coincide in $\cc H\Re$.
\end{cor}

\begin{proof}
Since $s$ is an $I$-weak equivalence, it follows that $rB^I$ is a
cylinder object for $rB$. The proof now follows from
Theorem~\ref{ccc} and Corollary~\ref{ddd}.
\end{proof}

\begin{exs}{\rm (1) Let $A\in\ring$. The group $GL_n(A)$ is defined as
$\kr(GL_n(\epsilon):GL_n(A^+)\to GL_n(\bb Z))$. Here $A^+=\bb Z\ps
A$ as a group and
   $$(n,a)(m,b)=(nm,nb+ma+ab).$$
We put $\epsilon:A^+\to\bb Z$ to be the augmentation
$\epsilon(n,a)=n$ and $GL(A):=\colim_n GL_n(A)$. The associated
functor $A\longmapsto GL(A)$ in $U\Re_{\bullet}$, pointed at the
unit element, denote by $\cc Gl$.

By definition, the Karoubi-Villamayor $K$-theory is defined as
   $$KV_n(A)=\pi_{n-1}(GL(A[\varDelta])),\ \ \ n\geq 1.$$
It follows from Theorem~\ref{ccc} that
   $$KV_n(A)=\Hom_{\Ho_{I,\bullet}(\Re)}((rA_+)\wedge S^{n-1},\cc Gl),\ \ \ n\geq 1.$$

(2) Let $\bb K^B(R)$ be a non-connective $K$-theory (simplicial)
$\Omega$-spectrum, functorial in $R$, where $R$ is a ring with
unit. We can extend $\bb K^B$ to all rings by the rule
   $$R\in\ring\longmapsto fibre(\bb K^B(R^+)\to\bb K^B(\bb Z)).$$
If $R$ has a unit this definition is consistent because then
$R^+\cong\bb Z\times R$.

The {\it homotopy $K$-theory\/} of $R\in\ring$ in the sense of
Weibel~\cite{W1} is given by the (fibrant) geometric realization
$KH(R)$ of the simplicial spectrum $\bb K^B(R[\varDelta])$. Note
that $KH(R)$ is an $\Omega$-spectrum. For $n\in\bb Z$, we shall
write $KH_n(R)$ for $\pi_nKH(R)$.

Let $K(A)$ denote the zeroth term of the spectrum $\bb K^B(A)$.
The corresponding functor $[A\in\Re\longmapsto K(A)]\in U\Re$
denote by $\cc K$. It is pointed at zero. It follows from
Theorem~\ref{ccc} that
   $$KH_n(A)=\Hom_{\Ho_{I,\bullet}(\Re)}((rA_+)\wedge S^{n},\cc K),\ \ \ n\geq 0.$$

}\end{exs}

\section{Homology theories on rings}\label{hom}

In this section we shall construct homology theories on rings.
Precisely, one will naturally associate to any pointed simplicial
functor $\cc X\in U\Re_\bullet$ a homology theory $\{H_n=H^{\cc
X,\ff F}_n\}_{n\geq 0}:\Re\to Sets$ depending on the family of
fibrations $\ff F$ of rings defined below. Such a homology theory is
defined by means of an explicitly constructed pointed simplicial
functor $Ex_{I,J}(\cc X)\in U\Re_\bullet$ and, by definition,
   $$H_n(A)=\pi_n(Ex_{I,J}(\cc X)(A))$$
for any $A\in\Re$ and $n\geq 0$. Moreover, there is a natural
transformation $\theta_{\cc X}:\cc X\to Ex_{I,J}(\cc X)$,
functorial in $\cc X$.

There is another formula for $H_n(A)$. A model category
$U\Re_{I,J,\bullet}$ is constructed and then
   $$H_n(A)=\Ho_{I,J,\bullet}(S^n\wedge rA,\cc X),$$
where $\Ho_{I,J,\bullet}$ stands for the homotopy category of
$U\Re_{I,J,\bullet}$.

Roughly speaking, we turn any pointed simplicial functor into a
homology theory. In this way the important simplicial functors
$\cc Gl$ and $\cc K$ give rise to the homology theories
$\{KV_n\mid\ff F=\textrm{$GL$-fibrations}\}$ and $\{KH_n\mid\ff
F=\textrm{surjective maps}\}$ respectively.

\subsection{Fibrations of rings}

\begin{defs}{\rm
Let $\Re$ be an admissible category of rings. A family $\ff F$ of
surjective homomorphisms of $\Re$ is called {\it fibrations\/} if
it meets the following axioms:

\begin{itemize}

\item[Ax 1)] for each $R$ in $\Re$, $R\to 0$ is in
$\ff F$;

\item[Ax 2)] $\ff F$ is closed under composition and any isomorphism is a fibration;

\item[Ax 3)] if the diagram
   $$\xymatrix{D\ar[r]^\rho\ar[d]_\sigma &A\ar[d]^f\\
               B\ar[r]^g &C}$$
is cartesian in $\Re$ and $g\in\ff F$, then $\rho\in\ff F$. Call
such squares {\it distinguished}. We also require that the
``degenerate square'' with only one entry, 0, in the upper
left-hand corner be a distinguished square;

\item[Ax 4)] any map $u$ in $\Re$ can be factored $u=pi$, where $p$ is a
fibration and $i$ is an $I$-weak equivalence.
\end{itemize}
Notice that the axioms imply that $\Re$ is closed under finite
direct products. We call a short exact sequence in $\Re$
   $$A\lra{g}B\lra{f}C$$
with $f\in\ff F$ a {\it $\ff F$-fibre sequence}.

$\ff F$ is said to be {\it saturated\/} if the homomorphism
$\partial_x^1:EA\to A$ is a fibration for any $A\in\Re$.

}\end{defs}

The trivial case is $\Re=\ff F=0$. A non-trivial example, $\Re\neq
0$, of fibrations is given by the surjective homomorphisms.
Indeed, the axioms Ax~1)-Ax~3) are trivial and Ax-4) follows from
Lemma~\ref{zzz} below.

Another important example of fibrations is defined by any left
exact functor. Recall that a functor $F:\ring\to Sets$ is {\it
left exact\/} if $F$ preserves finite limits. In particular, if
$A\to B\to C$ is a short exact sequence in $\ring$, then
   $$0\to FA\to FB\to FC$$
is an exact sequence of pointed sets (since the zero ring is a
zero object in $\ring$, it determines a unique element of $FA$).
Furthermore $F$ preserves cartesian squares.

For instance, any representable functor is left exact as well as
the functor (see Gersten~\cite{G})
   $$R\in\ring\longmapsto GL(R).$$

\begin{defs}{\rm
A surjective map $g:B\to C$ is said to be a {\it $F$-fibration\/}
(where $F:\ring\to Sets$ is a functor) if $F(E^n(g)):FE^nB\to FE^nC$
is surjective for all $n>0$. Observe that nothing is said about
$F(g):FB\to FC$. It follows that if the composite $fg$ of two maps
is a $F$-fibration, then so is $f$. If $F=GL$ we refer to
$F$-fibrations as {\it $GL$-fibrations}. We also note that the
family of all surjective homomorphisms is the family of
$F$-fibrations with $F$ sending a ring $A$ to itself.

}\end{defs}

\begin{lem}\label{zzz}
The collection of $F$-fibrations, where $F:\Re\to Sets$ is left
exact, enjoys the axioms Ax~1)-4) for fibrations on $\Re$ and is
saturated.
\end{lem}

\begin{proof}
The axioms Ax~1)-3) and the fact that $\ff F$ is saturated follow
from Gersten~\cite{G1}. Let us check Ax~4).

Let $u:A\to B$ be a homomorphism in $\Re$. Consider the following
commutative diagram
   $$\xymatrix{EB\ar@{ >->}[r]^{\nu}\ar@{=}[d] & A'\ar[d]^{\iota_1}\ar[r]^{\iota_2}& A\ar[d]^u\\
               EB \ar@{ >->}[r]^{\mu} & B[x]\ar[r]^{\partial^0_x}& B}$$
with $A'=A\times_BB[x]$. The map $i:A\to A'$,
$a\longmapsto(a,u(a))$, is split, $\iota_2i=1_A$, and obviously an
elementary homotopy equivalence. Hence it is an $I$-weak
equivalence.

Put $p:=\partial^1_x\circ\iota_1$. Then $p$ is surjective, because
any element $b\in B$ is the image of $(0,bx)$. By~\cite[2.3]{G1}
$F(E^n(p\nu))=F(E^n(\partial^1_x\mu))$, $n>0$, is a surjective
map. It follows that $F(E^n(p))$ is surjective. We see that $p$ is
a $F$-fibration.
\end{proof}

\subsection{The model category $U\Re_J$} We now introduce the class of excisive
functors on $\Re$. They look like flasque presheaves on a site
defined by a cd-structure in the sense of
Voevodsky~\cite[p.~14]{V}.

\begin{defs}{\rm
Let $\Re$ be an admissible category of rings and let $\ff F$ be a
family of fibrations. A simplicial functor $\cc X\in U\Re$ is
called {\it excisive\/} with respect to $\ff F$ if for any
distinguished square in $\Re$
   $$\xymatrix{
      D\ar[r]\ar[d]&A\ar[d]\\
      B\ar[r]&C
     }$$
the square of simplicial sets
   $$\xymatrix{
      \cc X(D)\ar[r]\ar[d]&{\cc X(A)}\ar[d]\\
      {\cc X(B)}\ar[r]&\cc X(C)
     }$$
is a homotopy pullback square. In the case of the degenerate
square the latter condition has to be understood in the sense that
$\cc X(0)$ is weakly equivalent to the homotopy pullback of the
empty diagram and is contractible. It immediately follows from the
definition that every pointed excisive object takes $\ff F$-fibre
sequences in $\Re$ to homotopy fibre sequences of simplicial sets.

}\end{defs}

\begin{exs}{\rm
Let $\ff F$ be the family of $GL$-fibrations. It follows
from~\cite{W2} that the simplicial functor
   $$A\in\Re\longmapsto Sing_*(\cc Gl)(A)=GL(A[\Delta])$$
is excisive.

The same is valid (see Weibel~\cite[Excision Theorem 2.2]{W1}) for
the homotopy $K$-theory simplicial functor
   $$A\in\Re\longmapsto Sing_*(\cc K)(A)$$
if $\ff F$ consists of all surjective homomorphisms.

}\end{exs}

Let $\alpha$ denote a distinguished square in $\Re$
   $$\xymatrix{
      D\ar[r]\ar[d]&A\ar[d]\\
      B\ar[r]&C
     }$$
and denote the pushout of the diagram
   $$\xymatrix{
      rC\ar[r]\ar[d]&rA\\
      rB}$$
by $P(\alpha)$. Notice that the obtained diagram is homotopy
pushout. There is a natural map $P(\alpha)\to rD$, and both
objects are cofibrant. In the case of the degenerate square this
map has to be understood as the map from the initial object
$\emptyset$ to $r0$.

We can localize $U\Re$ (respectively $U\Re_\bullet$) at the family
of maps
   $$J=\{P(\alpha)\to rD\mid\textrm{ $\alpha$ is a distinguished square}\}.$$
(respectively $J=\{P(\alpha)_+\to (rD)_+\}_\alpha$). The
corresponding $J$-localization will be denoted by $U\Re_J$
(respectively $U\Re_{J,\bullet}$). The weak equivalences (trivial
cofibrations) of $U\Re_{J}$ will be referred to as $J$-weak
equivalences ($J$-trivial cofibrations).

It follows that the square ``$r(\alpha)$''
   $$\xymatrix{
      rC\ar[r]\ar[d]&{rA}\ar[d]\\
      {rB}\ar[r]& rD}$$
with $\alpha$ a distinguished square is a homotopy pushout square
in $U\Re_J$.

The proof of the next lemma is straightforward.

\begin{lem}\label{zzo}
For any two rings $A,B\in\Re$, the natural map
   $$rA\sqcup rB\to r(A\times B)$$
is a $J$-weak equivalence. Therefore the simplicial set $\cc
X(A)\times\cc X(B)$ is weakly equivalent to the simplicial set
$\cc X(A\times B)$ for any $J$-local object $\cc X$. In
particular, the natural map
   $$rA\sqcup pt=rA\sqcup r0\to rA$$
is a $J$-weak equivalence.
\end{lem}

\begin{lem}\label{zz}
A simplicial functor $\cc X$ in $U\Re$ (respectively
$U\Re_\bullet$) is $J$-local \ifff it is fibrant and excisive.
\end{lem}

\begin{proof}
Straightforward.
\end{proof}

Now we define the mapping cylinder $\cyl(f)$ of a map $f:A\to B$
between cofibrant objects in a simplicial model category $\cc M$.
Let $A\otimes\varDelta^1$ denote the standard cylinder object for
$A$. One has a commutative diagram
   $$\xymatrix{A\sqcup A\ar[r]^(.55)\nabla\ar[d]_{i=i_0\sqcup i_1} &A\\
               A\otimes\varDelta^1\ar[ur]_\sigma}$$
in which $i$ is a cofibration and $\sigma$ is a weak
equivalence~\cite[9.5.14]{Hir}. Each $i_\epsilon$ must be a
trivial cofibration.

Form the pushout diagram
   $$\xymatrix{
      A\ar[r]^f\ar[d]_{i_0}&B\ar[d]^{i_{0*}}\\
      A\otimes\varDelta^1\ar[r]^{f_*} &\Cyl(f).
     }$$
Then $(f\sigma)\circ i_0=f$, and so there is a unique map
$q:\Cyl(f)\to B$ such that $qf_*=f\sigma$ and $qi_{0*}=1_B$. Put
$\cyl(f)=f_*i_1$; then $f=q\circ\cyl(f)$.

Since the objects $A, B, A\otimes\varDelta^1$ are cofibrant in
$\cc M$, it follows from~\cite[II.8.1]{GJ} that $\Cyl(f)$ is a
cofibrant object. Observe also that $q$ is a weak equivalence.

The map $\cyl(f)$ is a cofibration, since the diagram
   $$\xymatrix{
      A\sqcup A\ar[r]^{f\sqcup 1_A}\ar[d]_{i_0\sqcup i_1}&B\sqcup A\ar[d]^{i_{0*}\sqcup\cyl(f)}\\
      A\otimes\varDelta^1\ar[r]^{f_*} &\Cyl(f).
     }$$
is a pushout.

Given a distinguished square $\alpha$ let $P(\alpha)\to D_\alpha$
denote the cofibration $\cyl(P(\alpha)\to rD)$. We shall consider
the following set of maps
   $$\varLambda(J)=\{P(\alpha)\times\varDelta^n\bigsqcup_{P(\alpha)\times\partial\varDelta^n}
     D_\alpha\times\partial\varDelta^n\to D_\alpha\times\varDelta^n\}_{n\geq 0,\alpha}$$
In the pointed case one considers the set
   $$\varLambda(J)=\{P(\alpha)_+\wedge\varDelta^n_+\bigsqcup_{P(\alpha)_+\wedge\partial\varDelta_+^n}
     D_\alpha\wedge\partial\varDelta_+^n\to D_\alpha\wedge\varDelta_+^n\}_{n\geq 0,\alpha}$$
with $P(\alpha)_+\to D_\alpha$ the cofibration
$\cyl(P(\alpha)_+\to(rD)_+)$. It follows from~\cite[9.3.7(3)]{Hir}
that each map of $\varLambda(J)$ is a $J$-trivial cofibration. Let
$\cc C$ be a generating set of trivial cofibrations in $U\Re$ and
put $\varLambda:=\varLambda(J)\cup\cc C$.

\begin{prop}\label{zb}
A simplicial functor $\cc X$ in $U\Re$ (respectively
$U\Re_\bullet$) is $J$-local \ifff the map $\cc X\to *$ has the
right lifting property with respect to every element of
$\varLambda$.
\end{prop}

\begin{proof}
The proof is similar to~\cite[4.2.4]{Hir}. Use as
well~\cite[9.4.7]{Hir}.
\end{proof}

Observe that if an object $\cc X\in U\Re$ has the right lifting
property with respect to every element of $\varLambda(J)$ then it
is excisive (again use~\cite[9.4.7]{Hir}).

\subsection{The model category $U\Re_{I,J}$}

In this paragraph we shall construct the model category
$U\Re_{I,J}$. It is the localization of $U\Re$ with respect to the
maps from $I\cup J$. We start with definitions.

\begin{defs}{\rm
Let $\Re$ be an admissible category of rings and let $\ff F$ be a
family of fibrations. A simplicial functor $\cc X\in U\Re$ is
called {\it quasi-fibrant\/} with respect to $\ff F$ if it is
homotopy invariant and excisive.

Let $J$ be as above. The model category $U\Re_{I,J}$ is, by
definition, the Bousfield localization of $U\Re$ with respect to
$I\cup J$. The homotopy category of $U\Re_{I,J}$ will be denoted
by $\Ho_{I,J}(\Re)$. The weak equivalences (trivial cofibrations)
of $U\Re_{I,J}$ will be referred to as $(I,J)$-weak equivalences
($(I,J)$-trivial cofibrations).

An {\it $(I,J)$-resolution functor\/} is a pair
$(Ex_{I,J},\theta)$ consisting of a functor $Ex_{I,J}:U\Re\to
U\Re$ and a natural transformation $\theta:1\to Ex_{I,J}$ such
that for any $\cc X$ the object $Ex_{I,J}(\cc X)$ is quasi-fibrant
and the morphism $\cc X\to Ex_{I,J}(\cc X)$ is an $(I,J)$-trivial
cofibration.

}\end{defs}

\begin{lem}\label{z}
A simplicial functor $\cc X\in U\Re$ is $(I,J)$-local \ifff it is
fibrant, homotopy invariant and excisive.
\end{lem}

\begin{proof}
Straightforward.
\end{proof}

\subsubsection*{An ``explicit'' $(I,J)$-resolution functor} The
purpose of this paragraph is to construct an explicit
$(I,J)$-resolution functor. It is constructed inductive\-ly as
follows (cf. Morel-Voevodsky~\cite[p.~92]{MV}).

Given $\cc X\in U\Re$, let $\varLambda(J)$ be the set of
$J$-trivial cofibrations defined above and let $S$ be the set of
all commutative diagrams of the following form
   $$\xymatrix{P(\alpha)\times\varDelta^n\bigsqcup_{P(\alpha)\times\partial\varDelta^n}
     D_\alpha\times\partial\varDelta^n\ar[rr]\ar[d] &&Sing_*(\cc X)\ar[d]\\
     D_\alpha\times\varDelta^n\ar[rr]&& {*,}}$$
where $\alpha$ runs over distinguished squares. Construct a
pushout square
   $$\xymatrix{\bigsqcup_{S}[P(\alpha)\times\varDelta^n\bigsqcup_{P(\alpha)\times\partial\varDelta^n}
     D_\alpha\times\partial\varDelta^n]\ar[rr]\ar[d] &&Sing_*(\cc X)\ar[d]^{\xi_0}\\
     \bigsqcup_{S}D_\alpha\times\varDelta^n\ar[rr]&& \cc X_1.}$$
Since the left arrow is a $J$-trivial cofibration, then so is
$\xi_0$. We get a sequence of cofibrations
   $$\cc X\lra{\chi_0}Sing_*(\cc X)\lra{\xi_0}\cc X_1\lra{\chi_1}Sing_*(\cc X_1)$$
with $\chi_0,\chi_1$ $I$-trivial cofibrations. Repeating this
procedure, one obtains an infinite sequence of alternating
$I$-trivial cofibrations and $J$-trivial cofibrations
respectively,
   \begin{equation}\label{pp}
   \cdots Sing_*(\cc X_n)\lra{\xi_n}\cc X_{n+1}\xrightarrow{\chi_{n+1}}Sing_*(\cc X_{n+1})\cdots
   \end{equation}

\begin{prop}\label{nm}
Let $Ex_{I,J}(\cc X)$ denote a colimit of \eqref{pp} and let
$\theta_{\cc X}:\cc X\to Ex_{I,J}(\cc X)$ be the natural inclusion
which is functorial in $\cc X$. Then the pair $(Ex_{I,J},\theta)$
yields an $(I,J)$-resolution functor.
\end{prop}

\begin{proof}
The map $\theta_{\cc X}$ is a $(I,J)$-trivial cofibration
by~\cite[17.9.1]{Hir}. $Ex_{I,J}(\cc X)$ is plainly homotopy
invariant. To show that it is excisive, it is enough to check that
the map $Ex_{I,J}(\cc X)\to *$ has the right lifting property with
respect to all maps from $\varLambda(J)$. For this it suffices to
observe that both domains and codomains of maps in $\varLambda(J)$
commute with a colimit of~\eqref{pp}.
\end{proof}

An $(I,J)$-resolution functor $Ex_{I,J}(\cc X)$ with $\cc X\in
U\Re_\bullet$ is constructed in a similar way. The following
computes a fibrant replacement functor in $U\Re_{I,J}$
(respectively in $U\Re_{I,J,\bullet}$).

\begin{prop}\label{na}
Let $\vartheta_{\cc X}:\cc X\to R({\cc X})$ denote a fibrant
replacement functor in $U\Re$ (respectively in $U\Re_{\bullet}$).
Then the map $\vartheta_{Ex_{I,J}(\cc X)}\circ\theta_{\cc X}:\cc
X\longmapsto R(Ex_{I,J}(\cc X))$ yields a fibrant replacement
functor in $U\Re_{I,J}$ (respectively in $U\Re_{I,J,\bullet}$). That
is the object $R(Ex_{I,J}(\cc X))$ is $(I,J)$-local and the
composition
   $$\cc X\to {Ex_{I,J}(\cc X)}\to R(Ex_{I,J}(\cc X))$$
is an $(I,J)$-trivial cofibration.
\end{prop}

\begin{proof}
$R(Ex_{I,J}(\cc X))$ is plainly homotopy invariant. Given a
distinguished square $\alpha$, the square of simplicial sets
   $$\xymatrix{
      Ex_{I,J}(\cc X)(D)\ar[r]\ar[d]&{Ex_{I,J}(\cc X)(A)}\ar[d]\\
      {Ex_{I,J}(\cc X)(B)}\ar[r]&Ex_{I,J}(\cc X)(C)
     }$$
is a homotopy pullback square by Proposition~\ref{nm}. This square
is weakly equivalent to the square
   $$\xymatrix{
      R(Ex_{I,J}(\cc X))(D)\ar[r]\ar[d]&{R(Ex_{I,J}(\cc X))(A)}\ar[d]\\
      {R(Ex_{I,J}(\cc X))(B)}\ar[r]&R(Ex_{I,J}(\cc X))(C),
     }$$
and hence the latter square is a homotopy pullback square
by~\cite[13.3.13]{Hir}. Lemma~\ref{z} completes the proof.
\end{proof}

If we consider $rA$ as a pointed (at zero) simplicial functor then
the natural map $rA_+\to rA$ is a $J$-weak equivalence in
$U\Re_{\bullet}$ (see Lemma~\ref{zzo}). The proof of the following
statement is like that of Theorem~\ref{ccc}.

\begin{prop}\label{tan}
The natural map
   $$\pi_0(Ex_{I,J}(\cc X)(A))\to\Hom_{\Ho_{I,J}(\Re)}(rA,\cc X)$$
is a bijection for any $A\in\Re$ and $\cc X\in U\Re$. Moreover, if
$\cc X$ is pointed, then for any integer $n\geq 0$ and any
$A\in\Re$ the obvious map
   $$\pi_n(Ex_{I,J}(\cc X)(A))\to\Hom_{\Ho_{I,J,\bullet}(\Re)}(rA\wedge S^n,\cc X)$$
is a bijection, where $rA$ is supposed to be pointed at zero.
\end{prop}

\begin{cor}\label{an}
$(1)$ Suppose $\ff F$ consists of the $GL$-fibrations. Then for any
ring $A\in\Re$
   $$KV_n(A)=\Hom_{\Ho_{I,J,\bullet}(\Re)}(rA\wedge S^{n-1},\cc Gl),\ \ \ n\geq 1.$$

$(2)$ Suppose $\ff F$ consists of all surjective homomorphisms.
Then for any ring $A\in\Re$
   $$KH_n(A)=\Hom_{\Ho_{I,J,\bullet}(\Re)}(rA\wedge S^{n},\cc K),\ \ \ n\geq 0.$$
\end{cor}

\subsection{The Puppe sequence}
Throughout this paragraph the family of fibrations $\ff F$ is
supposed to be saturated. Let $g:B\to C$ be a ring homomorphism in
$\Re$. Consider the pullback of $g$ along the map
$\partial_x^1:EC=xC[x]\to C$,
   $$\xymatrix{P(g)\ar[r]^{g'}\ar[d]_{g_1} & EC\ar[d]^{\partial_x^1}\\
              B\ar[r]^g & C.}$$
Given a pointed quasi-fibrant simplicial functor $\cc X$, the
following lemma computes the homotopy type for $fibre(\cc
X(B)\to\cc X(C))$.

\begin{lem}\label{ccn}
If $\ff F$ is saturated and $\cc X$ is a pointed quasi-fibrant
simplicial functor, then the square of pointed simplicial sets
   $$\xymatrix{\cc X(P(g))\ar[r]\ar[d] & \cc X(EC)\simeq *\ar[d]\\
               \cc X(B)\ar[r] & \cc X(C)}$$
is homotopy pullback. In particular, it determines an exact
sequence of pointed sets at the middle point of the diagram
   $$[\cc X_0](P(g))\lra{g_{1*}}[\cc X_0](B)\lra{g_*}[\cc X_0](C).$$
\end{lem}

\begin{proof}
Easy.
\end{proof}

\begin{cor}\label{loop}
If $\ff F$ is saturated, $\cc X$ is a pointed quasi-fibrant
simplicial functor and $\Omega
A=(x^2-x)A[x]=\kr(EA\lra{\partial_x^1}A)$, $A\in\Re$, then $|\cc
X(\Omega A)|$ has the homotopy type of $\Omega|\cc X(A)|$. In
particular, $\pi_n(\cc X(A))=\pi_0(\cc X(\Omega^nA))$ for any
$n\geq 0$.
\end{cor}

Clearly one can iterate the construction of $P(g)$ to get the
diagram
   $$\xymatrix{\cdots\ar[r]&P(g_3)\ar@{=>}[r]^{g_4}\ar[d]_{g_3'}&P(g_2)\ar[r]^{g_2'}\ar@{=>}[d]_{g_3}&EP(g)\ar[d]^{g_3''}\\
               &EP(g_1)\ar[r]^{g_4''}&P(g_1)\ar@{=>}[r]^{g_2}\ar[d]_{g_1'}&P(g)\ar[r]^{g'}\ar@{=>}[d]_{g_1}&EC\ar[d]^{g_1''}\\
               &&EB\ar[r]^{g_2''}&B\ar@{=>}[r]^{g}&C.}$$
The latter diagram determines the {\it Puppe sequence of $g$}
   \begin{equation}\label{qu}
    \cdots\to P(g_n)\xrightarrow{g_{n+1}}P(g_{n-1})\lra{g_n}\cdots\to P(g)\lra{g_1}B\lra{g}C.
   \end{equation}
If we factor $g$ as $fi$ with $i$ a quasi-isomorphism and $f$ a
fibration, then using~\cite[II.9.10]{GJ} it is easy to show that
the Puppe sequence of $g$ is quasi-isomorphic to the Puppe
sequence of $f$.

\begin{prop}[cf. Gersten~\cite{G1}]\label{cb}
If $\ff F$ is saturated and $\cc X$ is a pointed quasi-fibrant
simplicial functor, then~\eqref{qu} gives rise to an exact
sequence of pointed sets
    \begin{gather*}\cdots\to [\cc X_0](P(g_n))\to[\cc X_0](P(g_{n-1}))\to\cdots\\
    \to [\cc X_0](P(g))\to[\cc X_0](B)\to[\cc X_0](C).
    \end{gather*}
\end{prop}

Gersten~\cite[2.9]{G1} constructs the same long exact sequence for
group valued left exact functors.

\subsection{Homology theories}

\begin{defs}[cf. Gersten~\cite{G}]{\rm
Let $\Re$ be an admissible category of rings and let $\ff F$ be a
family of fibrations. A {\it homology theory\/} $H_*$ on $\Re$
relative to $\ff F$ consists of
\begin{itemize}
\item[$(1)$] a family $\{H_n,n\geq 0\}$ of functors $H_n:\Re\to Sets_\bullet$
with $H_{n\geq 1}(A)$ a group,

\item[$(2)$] for every $\ff F$-fibre sequence
   $$A\lra f B\lra g C,$$
with $g\in\ff F$, morphisms
   $$H_{n+1}(C)\xrightarrow{\partial_{n+1}(g)}H_n(A),\ \ \ n\geq 0,$$
(we shall often write simply $\partial_{n+1}$ if $g$ is
understood) satisfying axioms \item[Ax 1)] $H_n(u)=H_n(v)$ for any
homotopic homomorphisms $u,v$ and any $n\geq 0$,

\item[Ax 2)] the morphism $\partial_{n+1}(g)$ of (2) is natural in
the sense that given a commutative diagram in $\Re$ with rows $\ff
F$-fibre sequences
    $$\xymatrix{A\ar[r]^f\ar[d]^a&B\ar[r]^g\ar[d]^b&C\ar[d]^c\\
                A'\ar[r]^{f'}&B'\ar[r]^{g'}&C'}$$
and with $g,g'\in\ff F$, then the diagram
    $$\xymatrix{H_{n+1}(C)\ar[rr]^{\partial_{n+1}(g)}\ar[d]_{H_{n+1}(c)}&&H_n(A)\ar[d]^{H_n(a)}\\
                H_{n+1}(C')\ar[rr]^{\partial_{n+1}(g')}&&H_n(A')}$$
is commutative for $n\geq 0$;

\item[Ax 3)] if  $A\lra f B\lra g C$ is an $\ff F$-fibre sequence
with $g\in\ff F$, then we have a long exact sequence of pointed
sets
   \begin{gather*}
    \cdots\to H_{n+1}(A)\xrightarrow{H_{n+1}(f)}H_{n+1}(B)\xrightarrow{H_{n+1}(g)}H_{n+1}(C)\\
    \xrightarrow{\partial_{n+1}(g)}H_n(A)\to\cdots\to H_0(B)\to
    H_0(C)\
   \end{gather*}
in the sense that the kernel (defined as the preimage of the
basepoint) is equal to the image at each spot.
\end{itemize}

}\end{defs}

We are now in a position to prove the following

\begin{thm}\label{ho}
To any pointed simplicial functor $\cc X$ on $\Re$ and any family
of fibrations $\ff F$ one naturally associates a homology theory.
It is defined as
   $$H_n(A):=\pi_n(Ex_{I,J}(\cc X)(A))=\Ho_{I,J,\bullet}(S^n\wedge rA,\cc X)$$
for any $A\in\Re$ and $n\geq 0$. Moreover, if $\ff F$ is saturated
then $H_n(A)=H_0(\Omega^nA)$. We also say that this homology
theory is {\em represented by $\cc X$}.
\end{thm}

\begin{proof}
It follows from Proposition~\ref{nm} that $Ex_{I,J}(\cc X)$ is a
quasi-fibrant object. Now our assertion easily follows from
Proposition~\ref{tan} and Corollary~\ref{loop}.
\end{proof}

It follows from Corollary~\ref{an} that the homology theories
associated to the functors $\cc Gl$ and $\cc K$ are the $KV$- and
$KH$-theories respectively.

\section{Derived categories of rings}\label{der}

In this section we introduce and study the left derived category
$D^-(\Re,\ff F)$ associated to any family of fibrations $\ff F$ on
$\Re$. It is obtained from the homotopy category $\cc H\Re$ by
inverting the quasi-isomorphisms introduced below. For this we
should first define a structure which is a bit weaker than the model
category structure on $\Re$ with respect to fibrations and
quasi-isomorphisms. Following Brown~\cite{B} this structure is
called the {\it category of fibrant objects}. It shares many
properties with model categories. If $\ff F$ is saturated (which is
always the case in practice), it follows from Theorem~\ref{rrr} that
$D^-(\Re,\ff F)$ is naturally left triangulated. The category of
left triangles meets the axioms which are versions for the axioms of
a triangulated category. The left triangulated structure as such is
a tool for producing homology theories on rings. The special case
when $\ff F$ consists of the surjective homomorphisms will be
discussed in section~\ref{kk}.

\subsection{Categories of fibrant objects}

\begin{defs}{\rm
I. Let $\cc A$ be a category with finite products and a final
object $e$. Assume that $\cc A$ has two distinguished classes of
maps, called {\it weak equivalences\/} and {\it fibrations}. A map
is called a {\it trivial fibration\/} if it is both a weak
equivalence and a fibration. We define a {\it path space\/} for an
object $B$ to be an object $B^I$ together with maps
   $$B\lra{s}B^I\xrightarrow{(d_0,d_1)}B\times B,$$
where $s$ is a weak equivalence, $(d_0,d_1)$ is a fibration, and
the composite is the diagonal map.

Following Brown~\cite{B}, we call $\cc A$ a {\it category of
fibrant objects\/} if the following axioms are satisfied.

(A) Let $f$ and $g$ be maps such that $gf$ is defined. If two of
$f$, $g$, $gf$ are weak equivalences then so is the third. Any
isomorphism is a weak equivalence.

(B) The composite of two fibrations is a fibration. Any
isomorphism is a fibration.

(C) Given a diagram
   $$A\bl u\longrightarrow C\bl v\longleftarrow B,$$
with $v$ a fibration (respectively a trivial fibration), the
pullback $A\times_CB$ exists and the map $A\times_CB\to A$ is a
fibration (respectively a trivial fibration).

(D) For any object $B$ in $\cc A$ there exists at least one path
space $B^I$ (not necessarily functorial in $B$).

(E) For any object $B$ the map $B\to e$ is a fibration.

}\end{defs}

Note that if the final object $e$ is also initial, then the
opposite category $\cc A^{\op}$ is a saturated Waldhausen category
(for precise definitions see~\cite{Wal,T}). The ``gluing axiom''
follows from~\cite[II.9.10]{GJ}. If $\cc A$ is small the
associated Waldhausen $K$-theory space of $\cc A^{\op}$ (see
Waldhausen~\cite{Wal}) will be denoted by $K\cc A$.

\begin{defs}{\rm
Let $\Re$ be an admissible category of rings and let $\ff F$ be a
family of fibrations. A homomorphism $A\to B$ in $\Re$ is said to be
a {\it $\ff F$-quasi-isomorphism\/} or just a {\it
quasi-isomorphism\/} if the map $rB\to rA$ is an $(I,J)$-weak
equivalence. This is equivalent to saying that the induced map
$H_*(A)\to H_*(B)$ is an isomorphism for every representable
homology theory $H_*$.

}\end{defs}

\begin{prop}\label{ttt}
Let $\Re$ be an admissible category of rings and let $\ff F$ be a
family of fibrations. Then it enjoys the axioms (A)-(E) for a
category of fibrant objects, where fibrations are the elements of
$\ff F$ and weak equivalences are quasi-isomorphisms.
\end{prop}

\begin{proof}
Clearly, the axioms (A), (B), (E) are satisfied. The axiom (D) is
a consequence of Ax~4). Indeed, let $B$ be a ring in $\Re$ and
consider homomorphisms
   $$B\lra{i}B[x]\xrightarrow{(\partial_x^0,\partial_x^1)}B\times B,$$
where $i$ is an $I$-weak equivalence and the composite is the
diagonal. By Ax~4) $(\partial_x^0,\partial_x^1)$ can be factored
$(\partial_x^0,\partial_x^1)=(d_0,d_1)\circ s'$, where $s'$ is an
$I$-weak equivalence and $(d_0,d_1)$ is a fibration. Put $s:=s'i$;
then the diagonal can be factored $diag=(d_0,d_1)\circ s$, hence
the axiom (D).

A pullback of a fibration is, by definition, a fibration. It
remains to check that a pullback of a trivial fibration is a
trivial fibration.

Suppose the square $\alpha$
   $$\xymatrix{D\ar[r]^\rho\ar[d]_\sigma &A\ar[d]^f\\
               B\ar[r]^g &C}$$
is distinguished in $\Re$ and $f$ is a trivial fibration. We must
show that $\sigma$ is a trivial fibration.

Since the morphism $r(f)$ is an $(I,J)$-trivial cofibration, then so
is the morphism $rB\to P(\alpha)=rA\bigsqcup_{rC}rB$
by~\cite[7.2.12]{Hir}. By definition, the morphism $P(\alpha)\to rD$
is a $J$-weak equivalence, whence our assertion follows.
\end{proof}

\begin{defs}{\rm
Let $\Re$ be an admissible category of rings and let $\ff F$ be a
family of fibrations. The {\it left derived category\/}
$D^-(\Re,\ff F)$ of $\Re$ with respect to $\ff F$ is the category
obtained from $\Re$ by inverting quasi-isomorphisms.

}\end{defs}

\begin{prop}\label{misha}
The family of quasi-isomorphisms in the category $\cc H\Re$ admits
a calculus of right fractions. The derived category $D^-(\Re,\ff
F)$ is obtained from $\cc H\Re$ by inverting the
quasi-isomorphisms.
\end{prop}

\begin{proof}
Let $C$ be a ring in $\Re$. By the proof of Proposition~\ref{ttt}
one can choose a path space $C^I$
   $$C\lra{s}C^I\xrightarrow{(d_0,d_1)}C\times C$$
with $s$ an $I$-weak equivalence. Consider a diagram
   \begin{equation}\label{mi}B\lra f C\bl t\longleftarrow A\end{equation}
with $t$ a quasi-isomorphism. Let $D$ be the limit of the diagram of
the solid arrows
   $$\xymatrix{&&D\ar@/_/@{.>}[dll]_{t'}\ar@{.>}[d]_h\ar@/^/@{.>}[drr]^{f'}\\
               B\ar[dr]^f&&C^I\ar[dl]_{d_0}\ar[dr]^{d_1}&&A\ar[dl]_t\\
               &C&&C}$$
It follows from~\cite[Lemma~3]{B} that $t'$ is a trivial fibration.
By Corollary~\ref{dd} $tf'=ft'$ in $\cc H\Re$. Thus \eqref{mi} fits into a
commutative diagram in $\cc H\Re$,
   $$\xymatrix{D\ar[r]^{f'}\ar[d]_{t'} & A\ar[d]^t\\
               B\ar[r]^f &C}$$
with $t'$ a quasi-isomorphism.

Given $f,g:A\rightrightarrows B$, suppose there is a
quasi-isomorphism $t:B\to C$ such that $tf=tg$ in $\cc H\Re$.
By~\cite[Propositions~1-2]{B} there is a quasi-isomorphism
$t':A'\to A$ such that $ft'$
is homotopic to $gt'$ by a homotopy $h:A'\to C^I$. It follows from
Corollary~\ref{dd} that $ft'=gt'$ in $\cc H\Re$, and hence $\cc
H\Re$ admits a calculus of right fractions.
\end{proof}

\begin{rem}{\rm
There is a generalization, due to Cisinski~\cite{C}, for the
notion of a category of fibrant objects: the ``cat\'egorie
d\'erivable \`a gauche''. For such a category Cisinski describes
(similar to Brown~\cite{B}) its derived category.  We can also
conform his construction to admissible categories of rings, but we
shall leave this to the interested reader.

}\end{rem}

\begin{question}{\rm
Let $X$ be an object of $D^-(\Re,\ff F)$. Is it true that the
functor
   $$[X,-]=\Hom_{D^-(\Re,\ff F)}(X,-)$$
is represented by $rX$? It is equivalent to the problem of whether
the functor $D^-(\Re,\ff F)\to\Ho_{I,J,\bullet}(\Re)$, induced by
$A\longmapsto rA$, is fully faithful.

}\end{question}

By Brown~\cite[Theorem~1]{B} there is a functor
$\Omega':D^-(\Re,\ff F)\to D^-(\Re,\ff F)$ such that for any ring
$B$ and any path space $B^I$, $\Omega' B$ can be canonically
identified with the fibre of $(d_0,d_1):B^I\to B\times B$.
Furthermore, $\Omega' B$ has a natural group structure. Let
$p:E\to B$ be a fibration with fibre $F$.
By~\cite[Proposition~3]{B} there is a natural map
$a:F\times\Omega' B\to F$ in $D^-(\Re,\ff F)$ which defines a
right action of the group $\Omega' B$ on $F$.

Following Quillen~\cite{Q}, we now define a {\it fibration
sequence\/} to be a diagram $F\to E\to B$ in $D^-(\Re,\ff F)$
together with an action $m:F\times\Omega' B\to F$ in $D^-(\Re,\ff
F)$ which are isomorphic to the diagram and action obtained from a
fibration in $\Re$. Let $A\in\Re$; the map
$m_*:[A,F]\times[A,\Omega' B]\to[A,F]$ will be denoted by
$(\alpha,\lambda)\longmapsto\alpha\cdot\lambda$.

\begin{thm}[Quillen~\cite{Q}, Brown~\cite{B}]\label{brown}
Given a fibration sequence
   $$F\lra i E\lra p B,\ \ \ F\times\Omega' B\to F,$$
there is an exact sequence in $D^-(\Re,\ff F)$
   $$\cdots\to \Omega' E\to\Omega' B\to F\to E\to B,$$
where exactness is interpreted as in~\cite[p.~I.3.8]{Q}. The
induced sequence
   $$\cdots\to[A,\Omega' E]\xrightarrow{(\Omega' p)_*}[A,\Omega' B]
     \lra{\partial_*}[A,F]\lra{i_*}[A,E]\lra{p_*}[A,B]$$
meets the following properties:

\begin{enumerate}
\item $(p_*)^{-1}(0)=\im i_*$;

\item $i_*\partial_*=0$ and
$i_*\alpha_1=i_*\alpha_2\Longleftrightarrow\alpha_2=\alpha_1\cdot\lambda$
for some $\lambda\in[A,\Omega' B]$;

\item $\partial_*(\Omega' p)_*=0$ and
$\partial_*\lambda_1=\partial_*\lambda_2\Longleftrightarrow
\lambda_2=(\Omega' p)_*\mu\lambda_1$ for some $\mu\in[A,\Omega' E]$
under the product in the group $[A,\Omega' B]$;

\item the sequence of group homomorphisms from $[A,\Omega' E]$ to
the left is exact in the usual sense.
\end{enumerate}
\end{thm}

\begin{cor}\label{on}
Let $\ff F$ be a saturated family of fibrations. Then $\Omega A$ is
canonically isomorphic to $\Omega' A$ for every $A\in\Re$. In
particular, $\Omega A$ is a group object in $D^-(\Re,\ff F)$.
\end{cor}

\begin{proof}
The proof is straightforward, using the preceding theorem and the
exact sequence $\Omega A\to EA\to A$ (recall that $EA$ is
contractible).
\end{proof}

Let $K(\Re,\ff F)$  denote the Waldhausen $K$-theory space
associated to a family of fibrations $\ff F$. Recall
from~\cite[p.~261]{T} that the group $K_0(\Re,\ff F)$ is abelian
and it is the free group on generators $[A]$ as $A$ runs over the
objects of $\Re$, modulo the two relations
\begin{itemize}
\item[$\diamond$] $[A]=[B]$ if there is a quasi-isomorphism
$A\lra{\sim}B$.

\item[$\diamond$] $[E]=[F]+[B]$ for all $\ff F$-fibre sequences
$F\rightarrowtail E\twoheadrightarrow B$.
\end{itemize}

The {\it Grothendieck group\/} $K_0(D^-(\Re,\ff F))$ of $D^-(\Re,\ff
F)$ is the free group on generators $(A)$ as $(A)$ runs over the
iso-classes of objects in $D^-(\Re,\ff F)$, modulo the relation:
$(E)=(F)+(B)$ for all fibration sequences $F\to E\to B$ in
$D^-(\Re,\ff F)$.

By~\cite[\S4, Proposition~4]{B} there is a fibration sequence
$\Omega' A\to 0\to A$ for any $A\in\Re$, hence $(\Omega' A)=-(A)$
in $K_0(D^-(\Re,\ff F))$. It follows that $(B)-(A)=(B\times\Omega'
A)$ and thus every element of $K_0(D^-(\Re,\ff F))$ is the class
$(A)$ of some $A$ in $\Re$. We leave to the reader to check that
the natural map
   $$K_0(\Re,\ff F)\to K_0(D^-(\Re,\ff F))$$
is an isomorphism of abelian groups.

\subsection{The left triangulated structure}
Fix a saturated family of fibrations $\ff F$ on $\Re$. In this
paragraph we define and study abstract properties of left triangles
in the derived category $D^-(\Re,\ff F)$.

The endofunctor $\Omega:\Re\to\Re$, $A\longmapsto\Omega
A=(x^2-x)A[x]$ respects quasi-isomorphisms. Indeed, let $f:A\to B$
be a quasi-isomorphism. Consider the following commutative
diagram:
   $$\xymatrix{\Omega A\ar@{ >->}[r]\ar[d]_{\Omega f}&EA\ar[d]_{E(f)}\ar@{->>}[r]^{\partial_x^1}&A\ar[d]^f\\
               \Omega B\ar@{ >->}[r]&EB\ar@{->>}[r]^{\partial_x^1}&B}$$
Since $EA,EB$ are isomorphic to zero in $D^-(\Re,\ff F)$, it
follows that $E(f)$ is a quasi-isomorphism. Then $\Omega f$ is a
quasi-isomorphism by~\cite[\S4, Lemma~3]{B}. Thus $\Omega$ can be
regarded as an endofunctor of $D^-(\Re,\ff F)$.

Given a fibration $g:A\to B$ with fibre $F$, consider the
commutative diagram as follows:
   $$\xymatrix{&\Omega B\ar@{=}[r]\ar@{ >->}[d]_j&\Omega B\ar@{ >->}[d]\\
               F\ar@{ >->}[r]^i\ar@{=}[d]&P(g)\ar@{->>}[r]\ar[d]^{g_1}&EB\ar[d]^{\partial_x^1}\\
               F\ar@{ >->}[r]^\iota&A\ar@{->>}[r]^g&B}$$
Since $EB$ is isomorphic to zero in $D^-(\Re,\ff F)$, it follows
from Theorem~\ref{brown} that $i$ is a quasi-isomorphism. We
deduce the sequence in $D^-(\Re,\ff F)$
   \begin{equation}\label{333}
    \Omega B\xrightarrow{i^{-1}\circ j}F\lra{\iota}A\lra{g}B.
   \end{equation}
We shall refer to such sequences as {\it standard left triangles}.
Any diagram in $D^-(\Re,\ff F)$ which is isomorphic to the latter
sequence will be called a {\it left triangle}. One must be careful
to note that $\Omega B'\to F'\lra{}A'\lra{}B'$ is isomorphic to a
standard triangle~\eqref{333} \ifff there is a commutative diagram
   $$\xymatrix{\Omega B\ar[r]\ar[d]_{\Omega b}&F\ar[d]_{f}\ar[r]&A\ar[d]^a\ar[r]&B\ar[d]^b\\
               \Omega B'\ar[r]&F'\ar[r]&A'\ar[r]&B'}$$
with $f,a,b$ isomorphisms in $D^-(\Re,\ff F)$.

It follows that the diagram
   $$\Omega B\lra{j}P(g)\lra{g_1}A\lra{g}B$$
is a left triangle. If $g$ is not a fibration then $g$ is factored
as $g=g'\ell$ with $g'$ a fibration and $\ell$ a
quasi-isomorphism. We get a commutative diagram
   $$\xymatrix{\Omega B\ar[r]\ar@{=}[d]&P(g)\ar[d]^t\ar[r]&A\ar[d]^\ell\ar[r]^g&B\ar@{=}[d]\\
               \Omega B\ar[r]&P(g')\ar[r]&A'\ar[r]^{g'}&B.}$$
The arrow $t$ is a quasi-isomorphism by~\cite[II.9.10]{GJ}. Hence
the upper sequence of the diagram is a left triangle. This also
verifies that any map in $D^-(\Re,\ff F)$ fits into a left
triangle.

For any ring $A$ the automorphism $\sigma=\sigma_A:\Omega A\to\Omega
A$ takes a polyno\-mi\-al $a(x)$ to $a(1-x)$. Notice that $\sigma$
is functorial in $A$ and $\sigma^2=1$. Given a morphism $\alpha$ in
$D^-(\Re,\ff F)$, by $-\Omega\alpha$ denote the morphism
$\Omega\alpha\circ\sigma=\sigma\circ\Omega\alpha$. For any $n\geq 1$
the morphism $(-1)^n\Omega\alpha$ means $\sigma^n\Omega\alpha$. Now
we want to check that for a standard left triangle
   $$\Omega B\xrightarrow{i^{-1}\circ j}F\lra{\iota}A\lra{g}B$$
the sequence
   $$\Omega A\xrightarrow{-\Omega g}\Omega B\xrightarrow{i^{-1}\circ j}F\lra{\iota}A$$
is a left triangle, too.

Consider the following diagram in $D^-(\Re,\ff F)$:
   $$\xymatrix{\Omega A\ar[r]^{-\Omega g}\ar@{=}[d]&\Omega B\ar[d]^\nu\ar[r]^{i^{-1}\circ j}&F\ar[d]^i\ar[r]^\iota&A\ar@{=}[d]\\
               \Omega A\ar[r]^\kappa&P(g_1)\ar[r]^{g_2}&P(g)\ar[r]^{g_1}&A,}$$
where $P(g_1)=P(g)\times_AEA$ and $\nu:\Omega B\to P(g_1)$ is the
natural inclusion taking $b(x)\in\Omega B$ to $((0,b(x)),0)$.
Moreover, $\nu$ is a quasi-isomorphism. The homomorphism $\kappa$
takes $a(x)\in\Omega A$ to $((0,0),a(x))\in P(g_1)$.

The right and the central squares are commutative. We want to
check that so is the left square. For this it is enough to show
that $\nu\circ\Omega g\circ\sigma$ is homotopic to $\kappa$. The
desired (elementary) homotopy is given by the homomorphism
   $$a(x)\in\Omega A\longmapsto((a(1-y),g(a(1-xy))),a(x(1-y)))\in P(g_1)[y].$$
It follows that the upper sequence is isomorphic to the lower
which is a left triangle by above.

Since every left triangle
   $$\Omega B'\lra{\gamma}F'\lra{\beta}A'\lra{\alpha}B'$$
is, by definition, isomorphic to a standard left triangle of the
form
   $$\Omega B\xrightarrow{i^{-1}\circ j}F\lra{\iota}A\lra{g}B,$$
we infer from above that the sequence
   $$\Omega A'\xrightarrow{-\Omega\alpha}\Omega B'\lra{\gamma}F'\lra{\beta}A'$$
is a left triangle.

Let $g:A\to B$ be a homomorphism in $\Re$ and let $g$ be factored
as $fi$ with $i:A\to A'$ a quasi-isomorphism and $f:A'\to B$ a
fibration. Then the Puppe sequence of $g$
    $$\cdots\to P(g_n)\xrightarrow{g_{n+1}}P(g_{n-1})\lra{g_n}\cdots\to P(g)\lra{g_1}A\lra{g}B$$
is quasi-isomorphic to the Puppe sequence of $f$. We also infer
from above that the latter is naturally isomorphic in $D^-(\Re,\ff
F)$ to the sequence
    $$\cdots\to\Omega F\xrightarrow{-\Omega\iota}\Omega A'\xrightarrow{-\Omega f}
    \Omega B\xrightarrow{j\circ i^{-1}} F\lra{\iota}A'\lra{f}B.$$
This isomorphism can be depicted as the following commutative
diagram in $D^-(\Re,\ff F)$ with the vertical arrows
quasi-isomorphisms (for simplicity we assume that $g$ is a
fibration). \footnotesize
   $$\xymatrix{\cdots\ar[r]&\Omega^2B\ar[r]^{-\Omega (i^{-1}j)}\ar[d]_{\Omega\nu}&\Omega F\ar[r]^{-\Omega\iota}\ar[d]_{\Omega i}
               &\Omega A\ar[r]^{-\Omega g}\ar@{=}[d]&\Omega B\ar[r]^{i^{-1}j}\ar[d]_{\nu}&F\ar[r]^{\iota}\ar[d]_i
               &A\ar[r]^g\ar@{=}[d]&B\ar@{=}[d]\\
               \cdots\ar[r]&\Omega P(g_1)\ar[r]^{-\Omega g_2}\ar@{=}[d]&\Omega P(g)\ar[r]^{-\Omega g_1}\ar@{=}[d]
               &\Omega A\ar[r]^{\kappa}\ar[d]_{i_1}&P(g_1)\ar[r]^{g_2}\ar@{=}[d]&P(g)\ar[r]^{g_1}\ar@{=}[d]
               &A\ar[r]^g\ar@{=}[d]&B\ar@{=}[d]\\
               \cdots\ar[r]&\Omega P(g_1)\ar[r]^{-\Omega g_2}\ar@{=}[d]&\Omega P(g)\ar[r]^{j_1}\ar[d]_{i_2}&P(g_2)\ar[r]^{g_3}\ar@{=}[d]
               &P(g_1)\ar[r]^{g_2}\ar@{=}[d]&P(g)\ar[r]^{g_1}\ar@{=}[d]
               &A\ar[r]^g\ar@{=}[d]&B\ar@{=}[d]\\
               \cdots\ar[r]&\Omega P(g_1)\ar[r]^{j_2}\ar[d]_{i_3}&P(g_3)\ar[r]^{g_4}\ar@{=}[d]&P(g_2)\ar[r]^{g_3}\ar@{=}[d]
               &P(g_1)\ar[r]^{g_2}\ar@{=}[d]&P(g)\ar[r]^{g_1}\ar@{=}[d]
               &A\ar[r]^g\ar@{=}[d]&B\ar@{=}[d]\\
               \cdots\ar[r]&P(g_4)\ar[d]\ar[r]^{g_5}&P(g_3)\ar[d]\ar[r]^{g_4}&P(g_2)\ar[d]\ar[r]^{g_3}&P(g_1)\ar[d]\ar[r]^{g_2}
               &P(g)\ar[d]\ar[r]^{g_1}&A\ar[d]\ar[r]^g&B\ar[d]\\
               &\vdots&\vdots&\vdots&\vdots&\vdots&\vdots&\vdots}$$
\normalsize

Now we are going to check the following property. Suppose we are
given two left triangles $\Omega
B\lra{\gamma}F\lra{\beta}A\lra{\alpha}B$ and $\Omega
B'\lra{\gamma'}F'\lra{\beta'}A'\lra{\alpha'}B'$ and two morphisms
$\phi:A\to A'$ and $\psi:B\to B'$ in $D^-(\Re,\ff F)$ with
$\psi\alpha=\alpha'\phi$. We claim that there exists a morphism
$\chi:F\to F'$ such that the triple $(\chi,\phi,\psi)$ is a morphism
from the first triangle to the second in the usual sense. It will
follow from the construction that $\chi$ is an isomorphism whenever
$\phi,\psi$ are.

Without loss of generality we can assume that the first left
triangle is the sequence $\Omega B\lra{j}P(g)\lra{g_1}A\lra{g}B$
and the second one is $\Omega
B'\lra{j'}P(g')\lra{g'_1}A'\lra{g'}B'$. Let $\psi=us^{-1}$ and
$\phi=vt^{-1}$.

Given two morphisms $g:A\to B$ and $s:Y\to B$ and a path space $B^I$ of $B$,
let $C:=Y\times_BB^I\times_BA$. One has a commutative diagram
   $$\xymatrix{&&C\ar@/_/[dll]_a\ar[d]_h\ar@/^/[drr]^c\\
               Y\ar[dr]^s&&B^I\ar[dl]_{d_0}\ar[dr]^{d_1}&&A\ar[dl]_g\\
               &B&&B}$$
The square
   $$\xymatrix{\ar @{}[dr] |{\textrm{h. comm.}}
               C\ar[d]_{a}\ar[r]^{c} & A \ar[d]^{g} \\
               Y\ar[r]_{s}&B}$$
is homotopy commutative by the homotopy $h$. Moreover, $c$ is a
quasi-isomorphism whenever $s$ is~\cite[\S2, Lemma~3]{B}. Maps from
a ring $D$ to $C$ correspond bijectively to data $(\ell,p,k)$ where
$\ell:D\to Y$ and $p:D\to A$ are maps and $k:D\to B^I$ is a homotopy
$gp\sim s\ell:D\to B$.

We can now construct the following commutative diagram in $D^-(\Re,\ff F)$:
   $$\xymatrix{&A\ar[rr]^g&&B\\
               X\ar[d]_v\ar[ur]^t&Z\ar[l]_p\ar[r]^q&C\ar[r]^a\ar[ul]_c&Y\ar[u]_s\ar[d]^u\\
               A'\ar[rrr]^{g'}&&&B',}$$
where $C:=Y\times_BB^I\times_BA$ and $Z:=X\times_AA^I\times_AC$. Let
$B'^I$ be a path space of $B'$. It follows from~\cite[\S2]{B} that
there exists a quasi-isomorphism $\ell:A''\to Z$ such that
$gcq\ell\sim saq\ell$ by a homotopy $k:A''\to B^I$ and $g'vp\ell\sim
uaq\ell$ by a homotopy $k':A''\to B'^I$.

We get the following commutative diagram in $\cc H\Re$:
   $$\xymatrix{A\ar[d]_g&\ar[l]_\tau A''\ar[r]^\alpha\ar[d]^\pi&A'\ar[d]^{g'}\\
               B&\ar[l]_s Y\ar[r]^u&B'}$$
with $\tau=cq\ell$ a quasi-isomorphism, $\pi=aq\ell$, and
$\alpha=vp\ell$.

\begin{lem}\label{rm}
Suppose we are given a homotopy commutative square with entries
$(X_0,Y,A_0,A_1)$
   $$\xymatrix{&& X_1\ar@{.>}[dl]^l\ar@/^/@{.>}[ddl]^{g''}\\
               \ar@{}[dr] |{\textrm{h. comm.}}
               X_0 \ar@/^/@{.>}[urr]^{f''}\ar[d]_{g'} \ar[r]^{f'} & Y \ar[d]^{g}\\
               A_0 \ar[r]_{f} & A_1}$$
and $gf'\sim fg'$ by a homotopy $k:X_0\to A_1^I$. Then there is
a $X_1$ and the dotted arrows $l,f'',g''$ such that the
square with entries $(X_0,X_1,A_0,A_1)$ is genuinely commutative,
$lf''=f'$, and $g''\sim gl$ by a homotopy $h:X_1\to A_1^I$. Moreover,
$l$ is a quasi-isomorphism.
\end{lem}

\begin{proof}
Let $X_1$ be the limit of the diagram of the solid arrows
   $$
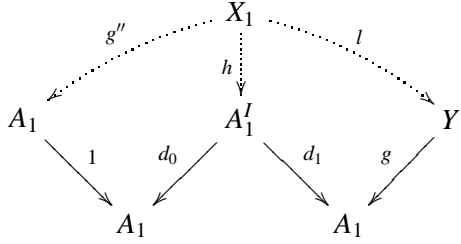
\xymatrix{&&X_1\ar@/_/@{.>}[dll]_{g''}\ar@{.>}[d]_h\ar@/^/@{.>}[drr]^l\\
               A_1\ar[dr]^1&&A_1^I\ar[dl]_{d_0}\ar[dr]^{d_1}&&Y\ar[dl]_g\\
               &A_1&&A_1}$$
The arrow $f''$ corresponds to the triple $(fg',f',k)$. Our assertion now
follows immediately.
\end{proof}

By Lemma~\ref{rm} one can construct a diagram in $\Re$
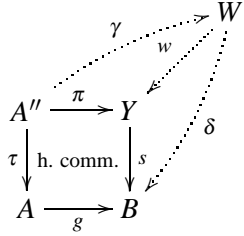
   $$\xymatrix{&& W\ar@{.>}[dl]_w\ar@/^/@{.>}[ddl]^{\delta}\\
               \ar@{}[dr] |{\textrm{h. comm.}}
               A''\ar@/^/@{.>}[urr]^{\gamma}\ar[d]_{\tau} \ar[r]^{\pi} & Y \ar[d]^{s}\\
               A \ar[r]_{g} & B}$$
resulting the diagram
   
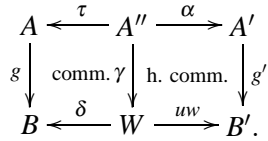
$$\xymatrix{\ar@{}[dr] |{\textrm{comm.}}
               A\ar[d]_g&\ar[l]_\tau A''\ar[r]^\alpha\ar[d]_\gamma\ar@{}[dr] |{\textrm{h. comm.}}&A'\ar[d]^{g'}\\
               B&\ar[l]_\delta W\ar[r]^{uw}&B'.}$$
In a similar way, one can construct a diagram
   $$
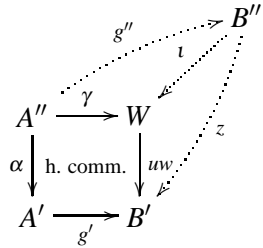
\xymatrix{&& B''\ar@{.>}[dl]_\iota\ar@/^/@{.>}[ddl]^{z}\\
               \ar@{}[dr] |{\textrm{h. comm.}}
               A''\ar@/^/@{.>}[urr]^{g''}\ar[d]_{\alpha} \ar[r]^{\gamma}&W\ar[d]^{uw}\\
               A'\ar[r]_{g'} & B'}$$
resulting the commutative diagram in $\Re$
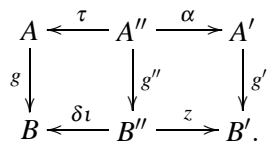
   $$\xymatrix{A\ar[d]_g&\ar[l]_\tau A''\ar[r]^\alpha\ar[d]^{g''}&A'\ar[d]^{g'}\\
               B&\ar[l]_{\delta\iota} B''\ar[r]^{z}&B'.}$$
Observe that $\phi=\alpha\tau^{-1}$ and $\psi=z(\delta\iota)^{-1}$. We
thus obtain the following commutative diagram in $\Re$
   $$\xymatrix{\Omega B\ar[r]^j&P(g)\ar[r]^{g_1}&A\ar[r]^g&B\\
               \ar[u]^{\Omega(\delta\iota)}\Omega B''\ar[d]_{\Omega z}\ar[r]^{j''}
               &\ar[u]P(g'')\ar[r]^{g_1''}\ar[d]&\ar[u]_\alpha A''\ar[r]^{g''}\ar[d]^\tau&\ar[u]_{\delta\iota}B''\ar[d]^z\\
               \Omega B'\ar[r]^{j'}&P(g')\ar[r]^{g_1'}&A'\ar[r]^{g'}&B'}$$
verifying the desired property.

We are now in a position to formulate the main result of the
paragraph.

\begin{thm}\label{rrr}
Let $\ff F$ be a saturated family of fibrations in $\Re$. Denote by
$\cc {L}tr(\Re,\ff F)$ the category of left triangles having the
usual set of morphisms from $\Omega C\lra{f}A\lra{g}B\lra{h}C$ to
$\Omega C'\lra{f'}A'\lra{g'}B'\lra{h'}C'$. Then $\cc {L}tr(\Re,\ff
F)$ is a {\it left triangulation\/} of $D^-(\Re,\ff F)$ in the sense
of Beligiannis-Marmaridis~\cite{BM}, i.e. it is closed under
isomorphisms and enjoys the following four axioms:

\begin{itemize}

\item[(LT1)] for any ring $A\in\Re$ the left triangle $0\lra{0}A\lra{1_A}A\lra{0}0$
belongs to $\cc {L}tr(\Re,\ff F)$ and for any morphism $h:B\to C$
there is a left triangle in $\cc {L}tr(\Re,\ff F)$ of the form
$\Omega C\lra{f}A\lra{g}B\lra{h}C$;

\item[(LT2)] for any left triangle $\Omega
C\lra{f}A\lra{g}B\lra{h}C$ in $\cc {L}tr(\Re,\ff F)$, the diagram
$\Omega B\xrightarrow{-\Omega h}\Omega C\lra{f}A\lra{g}B$ is also in
$\cc {L}tr(\Re,\ff F)$;

\item[(LT3)] for any two left triangles $\Omega
C\lra{f}A\lra{g}B\lra{h}C$, $\Omega
C'\lra{f'}A'\lra{g'}B'\lra{h'}C'$ in $\cc {L}tr(\Re,\ff F)$ and any
two morphisms $\beta:B\to B'$ and $\gamma:C\to C'$ of $D^-(\Re,\ff
F)$ with $\gamma h=h'\beta$, there is a morphism $\alpha:A\to A'$ of
$D^-(\Re,\ff F)$ such that the triple $(\alpha,\beta,\gamma)$ gives
a morphism from the first triangle to the second;

\item[(LT4)] any two morphisms $B\lra{h}C\lra{k}D$ of
$D^-(\Re,\ff F)$ can be fitted into a commutative diagram
   $$\xymatrix{&\Omega E\ar[d]^{f\circ\Omega\ell}\\
               \Omega C\ar[r]^f\ar[d]^{\Omega
               k}&A\ar[r]^g\ar[d]^\alpha&B\ar[r]^h\ar[d]_{1_B}&C\ar[d]^k\\
               \Omega D\ar[r]^j\ar[d]^{1_{\Omega
               D}}&F\ar[r]^m\ar[d]^{\beta}&B\ar[r]^{kh}\ar[d]_h&D\ar[d]^{1_D}\\
               \Omega D\ar[r]^i&E\ar[r]^\ell&C\ar[r]^k&D}$$
in which the rows and the second column from the left are left
triangles in $\cc Ltr(\Re,\ff F)$.
\end{itemize}
\end{thm}

The axiom (LT4) is a version of Verdier's octahedral axiom for left
triangles in $D^-(\Re,\ff F)$.

\begin{proof}
The axioms (LT1)-(LT3) are already checked above. Notice that the
morphism $\alpha$ in the axiom (LT3) is, by construction, an
isomorphism whenever $\beta,\gamma$ are. It remains to show (LT4).

Since every morphism in $D^-(\Re,\ff F)$ is of the form $p\circ
i\circ s^{-1}$ with $p$ a fibration and $i,s$ quasi-isomorphisms, it
follows that the composable morphisms $h,k$ fit into a commutative
diagram in $D^-(\Re,\ff F)$
   $$\xymatrix{B\ar[d]_{\cong}\ar[r]^h&C\ar[d]_{\cong}\ar[r]^k&D\ar[d]^{1_D}\\
               B'\ar[r]^p&C'\ar[r]^q&D}$$
with the vertical maps isomorphisms and $p,q$ fibrations in $\Re$.
It is routine to verify that (LT4) follows from the following fact
we are going to prove: any two fibrations $B\lra{h}C\lra{k}D$ of
$\Re$ can be fitted into a commutative diagram in $D^-(\Re,\ff F)$
   $$\xymatrix{&\Omega E\ar[d]^{f\circ\Omega\ell}\\
               \Omega C\ar[r]^f\ar[d]^{\Omega
               k}&A\ar[r]^g\ar[d]^\alpha&B\ar[r]^h\ar[d]_{1_B}&C\ar[d]^k\\
               \Omega D\ar[r]^v\ar[d]^{1_{\Omega
               D}}&F\ar[r]^m\ar[d]^{\beta}&B\ar[r]^{kh}\ar[d]_h&D\ar[d]^{1_D}\\
               \Omega D\ar[r]^u&E\ar[r]^\ell&C\ar[r]^k&D}$$
in which the rows are standard left triangles and the second column
from the left is a left triangle in $\cc Ltr(\Re,\ff F)$.

The horizontal standard triangles are constructed in a natural way
and then $\alpha,\beta$ exist by the universal property of
pullback diagrams. Note that $\beta$ is a fibration, because it is
base extension of the fibration $h$ along $\ell$. Moreover, the
sequence $A\lra\alpha F\lra\beta E$ is short exact. We have to
show that the sequence $\Omega
E\xrightarrow{f\circ\Omega\ell}A\lra\alpha F\lra\beta E$ is a left
triangle in $\cc Ltr(\Re,\ff F)$.

Recall that the map $f$ equals $i^{-1}\circ j$ with $i,j$ being
constructed as
   $$\xymatrix{&\Omega C\ar@{=}[r]\ar@{ >->}[d]_j&\Omega C\ar@{ >->}[d]\\
               A\ar@{ >->}[r]^i\ar@{=}[d]&P(h)\ar@{->>}[r]\ar[d]_{h_1}&EC\ar[d]^{\partial_x^1}\\
               A\ar@{ >->}[r]^g&B\ar@{->>}[r]^h&C.}$$
Let us construct a commutative diagram as follows
   $$\xymatrix{&\Omega E\ar@{=}[r]\ar@{ >->}[d]_\gamma&\Omega E\ar@{ >->}[d]\\
               A\ar@{ >->}[r]^\delta\ar@{=}[d]&P(\beta)\ar@{->>}[r]\ar[d]_{\beta_1}&E(E)\ar[d]^{\partial_x^1}\\
               A\ar@{ >->}[r]^\alpha&F\ar@{->>}[r]^\beta&E.}$$
It follows that $\delta$ is a quasi-isomorphism. Our assertion would
be proved if we show that the diagram
   $$\xymatrix{\Omega E\ar[rr]^{i^{-1}\circ j\circ\Omega\ell}\ar@{=}[d]&&A\ar[d]^{\delta}\ar[r]^\alpha&F\ar@{=}[d]\ar[r]^\beta&E\ar@{=}[d]\\
               \Omega E\ar[rr]^{\gamma}&&P(\beta)\ar[r]^{\beta_1}&F\ar[r]^\beta&E}$$
is commutative in $D^-(\Re,\ff F)$, because the lower sequence is a
left triangle. The left and central squares are commutative. It
remains to verify that $\delta\circ i^{-1}\circ
j\circ\Omega\ell=\gamma$.

By the universal property of pullback diagrams there exists a
homomorphism $\psi:P(\beta)\to P(h)$ making the diagram
   $$\xymatrix@!0{
     &E(E)\ar[rr]\ar'[d][dd] && E\ar[dd]^\ell\\
     P(\beta)\ar[ur]\ar[rr]\ar[dd]_{\psi} && F\ar[ur]\ar[dd]\\
     &EC\ar'[r][rr] && C\\
     P(h)\ar[rr]\ar[ur] && B\ar[ur]_h
     }$$
commutative. By construction, $\psi(f,e(x))=(m(f),\ell(e(x)))$ for
$(f,e(x))\in P(\beta)$. It follows that
$\psi\gamma=j\circ\Omega\ell$ and $\psi\delta=i$. Since $\delta,i$
are quasi-isomor\-phisms, then so is $\psi$. We have:
   $$\delta\circ i^{-1}\circ j\circ\Omega\ell=\delta i^{-1}\psi\gamma=\delta\delta^{-1}\gamma=\gamma.$$
Our theorem is proved.
\end{proof}

\begin{cor}\label{onn}
Let $\ff F$ be a saturated family of fibrations and $A\in\Re$. Then
the representable functor
   $$[A,-]=\Hom_{D^-(\Re,\ff F)}(A,-)$$
gives rise to a homology theory $H_*$ on $\Re$ with
$H_n(B)=[A,\Omega^nB]$, $n\geq 0$, and
   $$H_n(f)=
     \left\{
      \begin{array}{rcl}
       [A,(-1)^n\Omega (f)],\ n&\geq&1\\
       {[A,f]\hskip22pt},\ n&=&{0}
      \end{array}
      \right.$$
\end{cor}

\begin{proof}
By Corollary~\ref{on} $H_{n\geq 1}(B)$ is a group. Our assertion
would be proved we show that for any left triangle $\Omega
B\lra{f}F\lra{g}E\lra{h}B$ the induced sequence
   \begin{equation}\label{111}
    [A,\Omega B]\xrightarrow{\partial_1:=f_*}[A,F]\lra{g_*}[A,E]\lra{h_*}[A,B]
   \end{equation}
is an exact sequence of pointed sets. Since any left triangle is, by
definition, isomorphic to that induced by an $\ff F$-fibre sequence,
\eqref{111} is exact at the term $[A,E]$ by Theorem~\ref{brown}. By
(LT2) $\Omega E\xrightarrow{-\Omega h}\Omega B\lra{f}F\lra{g}E$ is a
left triangle. The same argument shows that \eqref{111} is exact at
$[A,F]$.
\end{proof}

\section{Stabilization}\label{stab}

Throughout this section $\Re$ is an admissible category of rings and
$\ff F$ is a saturated family of fibrations. There is a general
method of stabilizing the loop functor $\Omega$ (see
Heller~\cite{Hel}) and producing a triangulated category $D(\Re,\ff
F)$ from the left triangulated structure on $D^-(\Re,\ff F)$. We use
stabilization to define a $\bb Z$-graded bivariant homology theory
$k_*(A,B)$ on $\Re$, i.e. it is contravariant in the first variable
and covariant in the second and produces long exact sequences in
each variable out of $\ff F$-fibre sequences.

We start with preparations. First let us verify that
$\Omega^{n\geq 2}A$ are abelian group objects.

Let $B[x]\times_BB[x]:=\{(f(x),g(x))\mid f(1)=g(0)\}$ and let
$\wt\Omega B$ be the kernel of $(d_0,d_1):B[x]\times_BB[x]\to
B[x]$, $(f(x),g(x))\longmapsto (f(0),g(1))$. Denote by $\wt E$ the
fibred product of the diagram
   $$E\lra{\partial_x^1}B\bl{\partial_x^0}\longleftarrow B[x].$$
Since $\partial_x^1$ is a fibration then so is $pr_2:\wt E\to B[x]$.

\begin{lem}\label{zi}
The homomorphism $\alpha:\Omega B\to\wt\Omega B$, $f(x)\longmapsto
(f(x),0)$ is a quasi-isomorphism.
\end{lem}

\begin{proof}
Consider a commutative diagram in $\Re$
   $$\xymatrix{\Omega B\ar[r]\ar[d]^\alpha&\wt E\ar[r]^{pr_2}\ar[d]_{1}&B[x]\ar[d]^{\partial_x^1}\\
               \wt\Omega B\ar[r]\ar[d]^{\gamma}&\wt E\ar[r]^{p}\ar[d]_{pr_2}&B\ar[d]^{1}\\
               F\ar[r]^l&B[x]\ar[r]^{\partial_x^1}&B}$$
in which the rows are short exact. Note that $\gamma$ is a
fibration, because it is base extension of the fibration $pr_2$
along $l$. Thus the left column is a $\ff F$-fibre sequence.

The ring $F$ is quasi-isomorphic to 0, because it is isomorphic to
the contractible ring $E$. Therefore $\alpha$ is a
quasi-isomorphism by Theorem~\ref{brown}.
\end{proof}

Let us factorize $(\partial_x^0,\partial_x^1):B[x]\to B\times B$
as $B[x]\lra{i}B^I\lra{q}B\times B$, where $i$ is a $I$-weak
equivalence and $q$ is a fibration. Denote by $\Omega'B$ the fibre
of $(d_0,d_1):B^I\to B\times B$. The map $i$ induces a map
$u:\Omega B\to\Omega'B$. By Brown~\cite[p.~430-31]{B} one can
regard $\Omega'B$ as an object of $D^-(\Re,\ff F)$ well defined up
to canonical isomorphism. Moreover, we have a functor
$\Omega':\Re\to D^-(\Re,\ff F)$. This functor preservers
quasi-isomorphisms and so can be regarded as a functor
$\Omega':D^-(\Re,\ff F)\to D^-(\Re,\ff F)$.

Let $B^{2I}:=B^I\times_BB^I$ and let $\wt\Omega'B$ denote the
fibre of $(d_0,d_1):B^{2I}\to B\times B$. The map
$(i,i):B[x]\times_BB[x]\to B^{2I}$ yields a map $v:\wt\Omega
B\to\wt\Omega' B$. The maps of path spaces
$(1,sd_1),(sd_0,1):B^I\to B^{2I}$ induce two quasi-isomorphisms
$a,b:\Omega'B\to\wt\Omega'B$ taking $f\in\Omega'B$ to $(f,0)$ and
$(0,f)$ respectively. It follows from~\cite[\S4, Lemma~4]{B} that
$a=b$ in $D^-(\Re,\ff F)$.

\begin{lem}\label{zii}
The homomorphisms $u:\Omega B\to\Omega'B$ and $v:\wt\Omega
B\to\wt\Omega'B$ are quasi-isomorphisms.
\end{lem}

\begin{proof}
Consider a commutative diagram in $\Re$ with exact rows
$$\xymatrix{E\ar[r]\ar[d]^e&B[x]\ar[r]^{\partial_x^0}\ar[d]_i&B\ar[d]^{1}\\
               E'\ar[r]&B^I\ar[r]^{d_0}&B.}$$
Since $E,E'$ are quasi-isomorphic to zero, it follows that $e$ is a
quasi-isomorphism. Consider now a commutative diagram in $\Re$ with
exact rows
   $$\xymatrix{\Omega B\ar[r]\ar[d]^u&E\ar[r]^{\partial_x^1}\ar[d]&B\ar[d]^{1}\\
               \Omega'B\ar[r]&E'\ar[r]^{d_1}&B.}$$
By the proof of the factorization lemma in Brown~\cite{B} the map
$d_1$ is a fibration. It follows from~\cite[\S4, Lemma~3]{B} that
$u$ is a quasi-isomorphism. Since $v\alpha=au$ and $\alpha,a,u$
are quasi-isomorphisms (see Lemma~\ref{zi}), then so is $v$.
\end{proof}

Denote by $\beta:\Omega B\to\wt\Omega B$ the map taking
$f\in\Omega B$ to $(0,f)\in\wt\Omega B$. Since
$bu=v\beta,au=v\alpha$, $u$ and $v$ are quasi-isomorphisms by the
preceding lemma, and $a=b$ in $D^-(\Re,\ff F)$ we deduce the
following

\begin{cor}\label{zin}
$\alpha=\beta$ in $D^-(\Re,\ff F)$. In particular, $\beta$ is a
quasi-isomorphism.
\end{cor}

We now construct the following commutative diagram
   $$\xymatrix{\Omega B\times\Omega B\ar[r]^(.6)\omega\ar[d]_{(u,u)}&\wt\Omega B\ar[d]_v&\ar[l]_\alpha\Omega B\ar[d]^u\\
               \Omega'B\times\Omega'B\ar[r]^(.6)w&\wt\Omega'B&\ar[l]_a\Omega'B,}$$
where $\omega,w$ are obvious maps and the vertical maps are
quasi-isomorphisms by Lemma~\ref{zii}. Recall from Brown~\cite[p.
431]{B} that the map
   $$m_B:=a^{-1}w:\Omega'B\times\Omega'B\to\Omega'B$$
gives a group structure for $\Omega'B$ in $D^-(\Re,\ff F)$. The
map
   $$\mu_B:=\alpha^{-1}\omega:\Omega B\times\Omega B\to\Omega B$$
gives a group structure for $\Omega B$ in $D^-(\Re,\ff F)$,
because $\mu_B$ is isomorphic in $D^-(\Re,\ff F)$ to $m_B$ by
above.

\begin{lem}\label{ust}
For any ring $B\in\Re$ the homomorphism
$\tau:\Omega^2B\to\Omega^2B$, $\sum a_{ij}x^iy^j\mapsto\sum
a_{ij}x^jy^i$, is elementary homotopic to the identity.
\end{lem}

\begin{proof}
Any polynomial $f(x,y)\in\Omega^2 B$ can be written as
$f(x,y)=(x^2-x)(y^2-y)f'(x,y)$ for some (unique) polynomial
$f'(x,y)$. The desired elementary homotopy
$H:\Omega^2B\to\Omega^2B[t]$ is defined by
   $$(x^2-x)(y^2-y)f'(x,y)\bl H\longmapsto(x^2-x)(y^2-y)f'(tx+(1-t)y,(1-t)x+ty).$$
It follows that $d_0H=\tau$ and $d_1H=\id$.
\end{proof}

We are now in a position to prove the following

\begin{prop}\label{wii}
Let $\Re$ be an admissible category of rings and let $\ff F$ be a
saturated family of fibrations. Then for any ring $B\in\Re$ and
any $n\geq 2$ the ring $\Omega^{n}B$ is an abelian group object in
$D^-(\Re,\ff F)$.
\end{prop}

\begin{proof}
It will be sufficient to prove the claim for $\Omega^2B$. We use the
second coordinate to get the multiplication
$\Omega\mu_B:\Omega^2B\times\Omega^2B\to\Omega^2B$. First let us
show that $\Omega\mu_B=\mu_{\Omega B}$, i.e. the multiplications in
both coordinates agree.

The ring $\Omega\Omega'B$ is by construction consists of the
polynomials of the form $(x^2-x)\cdot[\sum_i(f_i(y),g_i(y))x^i]$
with each $(f_i(y),g_i(y))\in\Omega'B$. The ring $\Omega'\Omega B$
is by construction consists of the pairs $((y^2-y)\cdot[\sum_i
f_i(x)y^i],(y^2-y)\cdot[\sum_i g_i(x)y^i])$ with each
$(f_i(x),g_i(x))\in\Omega'B$.

Let $\tau'$ denote the homomorphism
   $$(x^2-x)\cdot[\sum_i(f_i(y),g_i(y))x^i]\longmapsto
     ((y^2-y)\cdot[\sum_i f_i(x)y^i],(y^2-y)\cdot[\sum_i g_i(x)y^i]).$$
Then the following diagram is commutative
   $$\xymatrix{\Omega^2B\times\Omega^2B\ar[d]_{\Omega\omega_B}\ar[rr]^{\tau\times\tau}&&\Omega^2B\times\Omega^2B\ar[d]^{\omega_{\Omega B}}\\
               \Omega\Omega'B\ar[rr]^{\tau'}&&\Omega'\Omega B\\
               \ar[u]^{\Omega\alpha_B}\Omega^2B\ar[rr]^{\tau}&&\ar[u]_{\alpha_{\Omega B}}\Omega^2B.}$$
By Lemma~\ref{ust} the upper and lower arrows equal identity in
$D^-(\Re,\ff F)$, hence $\Omega\mu_B=\mu_{\Omega B}$.

To verify that $\Omega^2B$ is an abelian group object in
$D^-(\Re,\ff F)$, one has to show that the diagram
   $$\xymatrix{\Omega^2B\times\Omega^2B\ar[rr]^T\ar[dr]_{\mu_{\Omega B}}&&\Omega^2B\times\Omega^2B\ar[dl]^{\mu_{\Omega B}}\\
               &\Omega^2B,}$$
in which $T$ is the isomorphism $(f,g)\longmapsto(g,f)$, is
commutative.

Let $T':\Omega\Omega'B\to\Omega'\Omega B$ denote the homomorphism
   $$(x^2-x)[\sum_i(f_i(y),g_i(y))x^i]\mapsto
     ((y^2-y)[\sum_i g_i(1-x)(1-y)^i],(y^2-y)[\sum_i f_i(1-x)(1-y)^i]).$$
Consider the diagram
   \begin{equation}\label{sta}
    \xymatrix{\Omega^2B\times\Omega^2B\ar[d]_{\Omega\omega_B}\ar[rr]^T&&\Omega^2B\times\Omega^2B\ar[d]^{\omega_{\Omega B}}\\
              \Omega\Omega'B\ar[rr]^{T'}&&\Omega'\Omega B\\
              \ar[u]^{\Omega\alpha_B}\Omega^2B\ar[rr]^{\id}&&\ar[u]_{\beta_{\Omega B}}\Omega^2B.}
   \end{equation}
We claim that is commutative in $D^-(\Re,\ff F)$.

Let $\sigma_x:\Omega^2B\to\Omega^2B$ (respectively $\sigma_y$) be
the homomorphism mapping $(x,y)$ to $(1-x,y)$ (respectively
$(x,1-y)$). We have
   $$\mu_{\Omega B}(\sigma_x\times\sigma_x)=\Omega\mu_B(\sigma_y\times\sigma_y),$$
and hence $\mu_{\Omega
B}=\Omega\mu_B(\sigma_x\sigma_y\times\sigma_x\sigma_y)=\mu_{\Omega
B}(\sigma_x\sigma_y\times\sigma_x\sigma_y)$. Then
   $$\beta^{-1}_{\Omega B}\omega_{\Omega B}T=\mu_{\Omega B}T=
     \mu_{\Omega B}T(\sigma_x\sigma_y\times\sigma_x\sigma_y)=\beta^{-1}_{\Omega B}T'\Omega\omega_B,$$
and so $\omega_{\Omega B}T=T'\Omega\omega_B$.

We also have
   $$T'\circ\Omega\alpha_B=T'\Omega\omega_B(1,0)^t=\omega_{\Omega
     B}T(1,0)^t=\omega_{\Omega B}(0,1)^t=\beta_{\Omega B}.$$
Here $(0,1)^t$ and $(1,0)^t$ denote the corresponding injections
$\Omega^2B\to\Omega^2B\times\Omega^2B$. The fact that $\beta_{\Omega
B}=\alpha_{\Omega B}$ (see Corollary~\ref{zin}) and that~\eqref{sta}
is commutative imply the desired abelian group structure on
$\Omega^2B$.
\end{proof}

\begin{cor}
Given two rings $A,B\in\Re$ and $n\geq 2$, the group ${D^-(\Re,\ff
F)}(A,\Omega^{n}(B))$ is abelian.
\end{cor}

We recall the construction of $D(\Re,\ff F)$ from
Heller~\cite{Hel}, which consists of formally inverting the
endofunctor $\Omega$. An object of $D(\Re,\ff F)$ is a pair
$(A,m)$ with $A\in D^-(\Re,\ff F)$ and $m\in\bb Z$. If $m,n\in\bb
Z$ then we consider the directed set $I_{m,n}=\{k\in\bb Z\mid
m,n\leq k\}$. The set of morphisms between $(A,m)$ and $(B,n)\in
D(\Re,\ff F)$ is defined by
   $$D(\Re,\ff F)[(A,m),(B,n)]:=\lp_{k\in I_{m,n}}D^-(\Re,\ff F)(\Omega^{k-m}(A),\Omega^{k-n}(B)).$$
Morphisms of $D(\Re,\ff F)$ are composed in the obvious fashion.
We define the {\it loop\/} automorphism on $D(\Re,\ff F)$ by
$\Omega(A,m):=(A,m-1)$. There is a natural functor $S:D^-(\Re,\ff
F)\to D(\Re,\ff F)$ defined by $A\longmapsto(A,0)$.

It follows from above that the category $D(\Re,\ff F)$ is
preadditive. Since it has finite direct products then it is
additive. We define a triangulation $\cc Tr(\Re,\ff F)$ of the
pair $(D(\Re,\ff F),\Omega)$ as follows. A sequence
   $$\Omega(A,l)\to (C,n)\to(B,m)\to(A,l)$$
belongs to $\cc Tr(\Re,\ff F)$ if there is an even integer $k$ and a
left triangle of representatives
$\Omega(\Omega^{k-l}(A))\to\Omega^{k-n}(C)\to\Omega^{k-m}(B)\to\Omega^{k-l}(A)$
in $D^-(\Re,\ff F)$. Clearly, the functor $S$ takes left triangles
in $D^-(\Re,\ff F)$ to triangles in $D(\Re,\ff F)$.

We are now in a position to prove the main result of the section.

\begin{thm}\label{r}
Let $\ff F$ be a saturated family of fibrations in $\Re$. Then
$\cc Tr(\Re,\ff F)$ is a triangulation of $D(\Re,\ff F)$ in the
classical sense of Verdier~\cite{Ver}.
\end{thm}

\begin{proof}
It is easy to see that $D(\Re,\ff F)$ is left triangulated, i.e.
$\cc Tr(\Re,\ff F)$ meets the axioms (LT1)-(LT4) of
Theorem~\ref{rrr}. By~\cite[p.~5]{BM} $D(\Re,\ff F)$ is
triangulated, because it is additive and the endofunctor $\Omega$
is invertible.
\end{proof}

We use the triangulated category $D(\Re,\ff F)$ to define a $\bb
Z$-graded bivariant homology theory depending on $(\Re,\ff F)$ as
follows:
   $$k_n(A,B)=k_n^{\Re,\ff F}(A,B):=D(\Re,\ff F)((A,0),(B,n)),\ \ \ n\in\bb Z.$$

\begin{cor}\label{rv}
For any $\ff F$-fibre sequence $A\to B\to C$ and any $D\in\Re$, we
have long exact sequences of abelian groups
   $$\cdots\to k_{n+1}(D,C)\to k_n(D,A)\to k_n(D,B)\to k_n(D,C)\to\cdots$$
and
   $$\cdots\to k_{n+1}(A,D)\to k_n(C,D)\to k_n(B,D)\to k_n(A,D)\to\cdots$$
\end{cor}

\section{The triangulated category $kk$}\label{kk}

Motivated by ideas and work of J. Cuntz on bivariant $K$-theory of
locally convex algebras (see~\cite{Cu,Cu1}), Corti\~nas and
Thom~\cite{CT} construct a bivariant homology theory $kk_*(A,B)$ for
algebras over a unital ground ring $H$. It is Morita invariant,
homotopy invariant, excisive $K$-theory of algebras, which is
universal in the sense that it maps uniquely to any other such
theory. This bivariant $K$-theory is defined similar to the
bivariant homology theory $k_*(A,B)$ discussed in the previous
section. Namely, a triangulated category $kk$ whose objects are the
$H$-algebras without unit is constructed and then set
$kk_n(A,B)=kk(A,\Omega^{n}B)$, $n\in\bb Z$.

We make use of our machinery developed in the preceding sections to
study various triangulated structures on admissible categories of
rings which are not necessarily small. As an application, we give
another description of the triangulated category $kk$. Throughout
this section the class $\ff F$ of fibrations consists of the
surjective homomorphisms.

Let $\Re$ be an arbitrary not necessarily small admissible category
of rings and let $\ff W$ be any subcategory of homomorphisms
containing the $I$-weak equivalences such that the triple $(\Re,\ff
W, \ff F)$ is a Brown category. Let $D^-(\Re,\ff W)$ be the category
obtained from $\Re$ by inverting the weak equivalences. Then
$\Omega:\Re\to\Re$ yields a loop functor on $D^-(\Re,\ff W)$. Let us
define the category of left triangles $\cc {L}tr(\Re,\ff W)$ similar
to $\cc {L}tr(\Re,\ff F)$. Then the following is true.

\begin{thm}\label{rrr1}
$\cc {L}tr(\Re,\ff W)$ determines a left triangulation of
$D^-(\Re,\ff W)$. The stabilization procedure of the loop functor
$\Omega$ described in the previous section yields a triangulated
category $D(\Re,\ff W)$ whose objects and morphisms are defined
similar to those of $D(\Re,\ff F)$.
\end{thm}

\begin{proof}
The proof repeats those of Theorems~\ref{rrr} and~\ref{r} word for
word if we replace in appropriate places path spaces $B^I$ with the
functorial path space $B[x]$ for $B\in(\Re,\ff W, \ff F)$.
\end{proof}

\begin{rem}{\rm
Theorem~\ref{rrr1} says that construction of $D(\Re,\ff F)$ is
formal and can be defined in a more general setting whenever $\ff F$
consists of the surjective homomorphisms.

}\end{rem}

\begin{cor}\label{rv1}
$k_*(A,B):=k_*^{\Re,\ff W}(A,B)=D(\Re,\ff W)(A,\Omega^*B)$
determines a bivariant homology theory on $\Re$. Moreover, for any
short exact sequence $A\to B\to C$ and any $D\in\Re$, we have long
exact sequences of abelian groups
   $$\cdots\to k_{n+1}(D,C)\to k_n(D,A)\to k_n(D,B)\to k_n(D,C)\to\cdots$$
and
   $$\cdots\to k_{n+1}(A,D)\to k_n(C,D)\to k_n(B,D)\to k_n(A,D)\to\cdots$$
\end{cor}

We consider associative, not necessarily unital or central algebras
over a fixed unital, not necessarily commutative ring $H$; we write
$\aha$ for the category of such algebras. By forgetting structure,
we can embed $\aha$ faithfully into each of the categories of
$H$-bimodules, abelian groups and sets. Fix one of these underlying
categories, call it $\cc U$, and let $F:\aha\to\cc U$ be the
forgetful functor. Let $\cc E$ be the class of all exact sequences
of $H$-algebras
\begin{equation}\label{extensionE}
(E):0\to A\to B\to C\to 0
\end{equation}
such that $F(B)\to F(C)$ is a split surjection.

\begin{defs}[Corti\~nas-Thom~\cite{CT}]{\rm
Given a triangulated category $(\cc T,\Omega)$, an {\it excisive
homology theory\/} on $\aha$ with values in $\cc T$ consists of a
functor $X:\aha\to \cc T$, together with a collection
$\{\partial_E:E\in\cc E\}$ of maps $\partial_E^X=\partial_E\in{\cc
T}(\Omega X(C), X(A))$. The maps $\partial_E$ are to satisfy the
following requirements. \sn \noindent{i)} For all $E\in \cc E$ as
above,
\[
\xymatrix{\Omega
X(C)\ar[r]^{\partial_E}&X(A)\ar[r]^{X(f)}&X(B)\ar[r]^{X(g)}& X(C)}
\]
is a distinguished triangle in $\cc T$. \sn \noindent{ii)} If
\[
\xymatrix{(E): &A\ar[r]^f\ar[d]_\alpha& B\ar[r]^g\ar[d]_\beta& C\ar[d]_\gamma\\
          (E'):&A'\ar[r]^{f'}& B'\ar[r]^{g'}& C'}
\]
is a map of extensions, then the following diagram commutes
\[
\xymatrix{{\Omega} X(C)\ar[d]_{{\Omega} X(\gamma)}\ar[r]^{\partial_E}& X(A)\ar[d]^{X(\alpha)}\\
 \Omega X(C')\ar[r]_{\partial_{E'}}& X(A).}
\]

Let $\iota_\infty:A\to M_\infty A$ be the natural inclusion from $A$
to $M_\infty A=\bigcup_nM_nA$, the union of matrix rings. An
excisive, homotopy invariant homology theory $X:\aha\to\cc T$ is
{\it $M_\infty$-stable\/} if for every $A\in\aha$, it maps the
inclusion $\iota_\infty:A\to M_\infty A$ to an isomorphism. Note
that if $X$ is $M_\infty$-stable, and $n\geq 1$, then $X$ maps the
inclusion $\iota_n:A\to M_nA$ to an isomorphism.

}\end{defs}

The homotopy invariant, $M_\infty$-stable, excisive homology
theories form a category, where a homomorphism between the theories
$X:\aha\to\cc T$ and $Y:\aha\to\cc S$ is a triangulated functor
$G:\cc T\to\cc S$ such that
   $$\xymatrix{{\aha}\ar[dr]_Y\ar[r]^X&{\cc T}\ar[d]^G\\
                                  &{\cc S}}$$
commutes, and such that for every extension \eqref{extensionE}, the
natural isomorphism $\phi:G(\Omega X(C))\to \Omega Y(C)$ makes the
following into a commutative diagram
\begin{equation}\label{phi-partial}
\xymatrix{G(\Omega X(C))\ar[r]^(.6){G(\partial^X_E)}\ar[d]_\phi&Y(A)\\
             \Omega Y(C).\ar[ur]_{\partial^Y_E}& }
\end{equation}

\begin{thm}[Corti\~nas-Thom~\cite{CT}]\label{untr}
The category of homotopy invariant, $M_\infty$-stable, excisive
homology theories has an initial object $\ell:\aha\to kk$. The
triangulated category $kk$ has the same objects and the same
endofunctor $\Omega$ as $\aha$. Furthermore,
   $$kk_*(A,B)=kk(A,\Omega^*B)$$
gives rise to a $M_\infty$-stable, homotopy invariant, excisive,
bivariant homology theory of algebras.
\end{thm}

The preceding theorem is used to define a natural map
   $$kk_*(A,B)\to KH_*(A,B),$$
where $KH_*(A,B)$ is the bivariant theory generated by the homotopy
$K$-theory $KH$. A result of Corti\~nas-Thom~\cite{CT} states that
this map is an isomorphism when $A=H$ is commutative and $B$ is a
central $H$-algebra, i.e. $kk_*(H,B)=KH_*(B)$. When $H$ is a field
of characteristic zero and $A,B$ are central $H$-algebras, they also
obtain in this way a product preserving Chern character to bivariant
periodic cyclic cohomology
   $$ch_*:kk_*(A,B)\to HP^*(A,B).$$

Let $\ff W_{CT}$ be the class of homomorphisms $f$ in $\aha$ such
that $X(f)$ is an isomorphism for any homotopy invariant,
$M_\infty$-stable, excisive homology theory $X:\aha\to\cc T$. It is
directly verified that the triple $(\aha,\ff W_{CT},\ff F)$ meets
the axioms for a Brown category.

We are now in a position to prove the main result of this section.

\begin{thm}[Comparison]\label{unr}
There is a natural triangulated equivalence of the triangulated
categories $kk$ and $D(\aha,\ff W_{CT})$.
\end{thm}

\begin{proof}
Let $\alpha:\Re\to D^-(\Re,\ff W),S:D^-(\Re,\ff W)\to D(\Re,\ff W)$
be the canonical functors and let $E$ be the extension
\eqref{extensionE}. Define $\partial_E\in D(\Re,\ff W)(\Omega C,A)$
as the class of the canonically defined morphism~\eqref{333}
$i^{-1}\circ j\in D^-(\Re,\ff W)(\Omega C,A)$. Then
$\iota:=S\alpha:\aha\to D(\aha,\ff W_{CT})$, together with
$\{\partial_E\}_{E\in\cc E}$ is a homotopy invariant,
$M_\infty$-stable, excisive homology theory. By Theorem~\ref{untr}
there is a unique morphism $G:kk\to D(\Re,\ff W)$ of homology
theories such that $G\ell=\iota$. We claim that $G$ is an
equivalence of categories.

Since $\ell$ takes weak equivalences to isomorphisms, there is a
unique functor $F:D^-(\aha,\ff W_{CT})\to kk$ such that
$F\circ\alpha=\ell$. We have $S\alpha=G\ell=GF\alpha$. It follows
that $S=GF$, and hence $G$ is full.

By~\cite[1.1]{Hel} $F$ is uniquely extended to a functor
$H:D(\aha,\ff W_{CT})\to\cc T$ such that $H\circ S=F$.
By~\cite[6.5.1]{CT} a diagram
   \[\xymatrix{{\Omega} C \ar[r]& A \ar[r] & B \ar[r] & C }\]
of morphisms in $kk$ is a distinguished triangle if it is isomorphic
in $kk$ to the path sequence
   \[\xymatrix{\Omega B'\ar[r]^{\ell(j)}&P(f)\ar[r]^{\ell(\pi_f)}&A'\ar[r]^{\ell(f)}&B'}\]
associated with a homomorphism $f:A'\to B'$ of $H$-algebras. We see
that $F$ takes left triangles in $D^-(\aha,\ff W_{CT})$ to triangles
in $kk$. Therefore $H$ takes triangles in $D(\aha,\ff W_{CT})$ to
triangles. The same argument as in the proof of~\cite[6.6.2]{CT}
shows that $H$ must be a morphism of homotopy invariant,
$M_\infty$-stable, excisive homology theories. By uniqueness
$HG=\id_{kk}$, and so $G$ is faithful, as was to be proved.
\end{proof}

\begin{cor}\label{unrrr}
$\ff W_{CT}$ is the smallest class of weak equivalences containing
the homomorphisms of $H$-algebras $A\to A[x],\iota_\infty:A\to
M_\infty A$, that is $\ff W_{CT}\subseteq\ff W$ with $\ff W$ being
any class of weak equivalences containing $A\to
A[x],\iota_\infty:A\to M_\infty A$ such that the triple $(\aha,\ff
W,\ff F)$ is a Brown category.
\end{cor}

\begin{proof}
Let $\ff W$ be the smallest class of weak equivalences containing
$A\to A[x],\iota_\infty:A\to M_\infty A$ such that the triple
$(\aha,\ff W,\ff F)$ is a Brown category. Then $\ff W\subseteq\ff
W_{CT}$. By Theorem~\ref{rrr1} the canonical functor $\aha\to
D(\aha,\ff W)$ yields a homotopy invariant, $M_\infty$-stable,
excisive homology theory. Therefore $\ff W_{CT}\subseteq\ff W$.
\end{proof}

We infer from the preceding theorem that $kk$ does not depend on the
choice of the underlying category $\cc U$. For further properties of
the category $kk$ we refer the reader to~\cite{CT}.

\section{Addendum}

When $\cc M$ is a model category and $S$ a set of maps between
cofibrant objects, we shall produce a new model structure on $\cc
M$ in which the maps $S$ are weak equivalences. The new model
structure is called the {\it Bousfield localization\/} or just
localization of the old one. A theorem of Hirschhorn says that
when $\cc M$ is a ``sufficiently nice'' model category one can
localize at any set of maps. ``Sufficiently nice'' entails being
cofibrantly generated together with having certain other
finiteness properties; the exact notion is that of a {\it
cellular\/} model category. We do not recall the definition here,
but refer the reader to~\cite{Hir}. Suffice it to say that all the
model categories we encounter in this paper are cellular.

For simplicity we shall from now on assume that all model
categories are simplicial. This is not strictly necessary, but it
allows us to avoid a certain machinery required for dealing with
the general case (see~\cite{Hir} for details).

Since all model categories we shall consider are simplicial, we do
not make use of the homotopy function complex $\map(X,Y)$ defined
in~\cite{Hir}. Indeed, let $\cc M$ be a simplicial model category
with simplicial mapping object $\Map$, and let $X$ and $Y$ be two
objects of $\cc M$. If $\wt X\to X$ is a cofibrant replacement of
$X$ and $Y\to\wh Y$ is a fibrant replacement of $Y$, then
$\map(X,Y)$ is homotopy equivalent to $\Map(\wt X,\wh Y)$.
Consequently, one can recast the localization theory of $\cc M$ in
terms of the simplicial mapping object instead of the homotopy
function complex.

\begin{defs}{\rm
Let $\cc M$ be a simplicial model category and let $S$ be a set of
maps between cofibrant objects.

\begin{enumerate}

\item An {\it $S$-local object\/} of $\cc M$ is a fibrant object
$X$ such that for every map $A \to B$ in $S$, the induced map of
$\Map(B,X)\to\Map(A,X)$ is a weak equivalence of simplicial sets.

\item An {\it $S$-local equivalence\/} is a map $A\to B$ such
that $\Map(B,X) \to\Map(A,X)$ is a weak equivalence for every
$S$-local object $X$.

\end{enumerate}
}\end{defs}

In words, the $S$-local objects are the ones which see every map
in $S$ as if it were a weak equivalence.  The $S$-local
equivalences are those maps which are seen as weak equivalences by
every $S$-local object.

\begin{thm}[Hirschhorn~\cite{Hir}]
Let $\cc M$ be a cellular, simplicial model category and let $S$
be a set of maps between cofibrant objects. Then there exists a
new model structure on $\cc M$ in which

\begin{enumerate}

\item the weak equivalences are the $S$-local equivalences;

\item the cofibrations in $\cc M/S$ are the same as those in $\cc
M$;

\item the fibrations are the maps having the
right-lifting-property with respect to cofibrations which are also
$S$-local equivalences.
\end{enumerate}
Left Quillen functors from $\cc M/S$ to $\cc D$ are in one to one
correspondence with left Quillen functors $\varPhi:\cc M\to\cc D$
such that $\varPhi(f)$ is a weak equivalence for all $f\in S$. In
addition, the fibrant objects of $\cc M$ are precisely the
$S$-local objects, and this new model structure is again cellular
and simplicial.

\end{thm}

The model category whose existence is guaranteed by the above
theorem is called {\it $S$-localization\/} of $\cc M$. The
underlying category is the same as that of $\cc M$, but there are
more trivial cofibrations (and hence fewer fibrations). We
sometimes use $\cc M/S$ to denote the $S$-localization.

Note that the identity maps yield a Quillen pair $\cc M
\rightleftarrows\cc M/S$, where the left Quillen functor is the
map $\id:\cc M\to\cc M/S$.


\begin{thebibliography}{99}

\bibitem{Bass}{\it H. Bass}, Algebraic K-theory, Benjamin, New
         York-Amsterdam, 1968, xx+762 p.
\bibitem{BM} {\it A. Beligiannis, N. Marmaridis}, Left triangulated
         categories arising from contravariantly finite subcategories,  Comm.
         Algebra 22(12) (1994), 5021-5036.
\bibitem{B}{\it K. S. Brown}, Abstract homotopy theory and generalized
            sheaf cohomology, Trans. Amer. Math. Soc. 186~(1973),
            419-458.
\bibitem{C}{\it D.-C. Cisinski}, Cat\'egories d\'erivables,
         preprint, 2002 (available at {\tt
         www-math.univ-paris13.fr/$\sim$cisinski}).
\bibitem{CT}{\it G. Corti\~{n}as, A. Thom}, Bivariant algebraic
         K-theory, preprint math.KT/0603531.
\bibitem{Cu} {\it J. Cuntz}, Bivariant K-theory and the Weyl algebra,
         K-theory 35 (2005), 93-137.
\bibitem{Cu1} {\it J. Cuntz, A. Thom}, Algebraic K-theory and locally convex algebras, Math. Ann. 334 (2006), 339-371.
\bibitem{D}{\it D. Dugger}, Sheaves and homotopy theory, preprint,
         1999 (available at {\tt darkwing.uoregon.edu/$\sim$ddugger}).
\bibitem{G}{\it S. M. Gersten}, On Mayer-Vietoris functors and algebraic K-theory,
         J. Algebra 18 (1971), 51-88.
\bibitem{G1}{\it S. M. Gersten}, Homotopy theory of rings,
         J. Algebra 19 (1971), 396-415.
\bibitem{GJ}{\it P. G. Goerss, J. F. Jardine}, Simplicial homotopy
        theory, Progress in Mathematics 174, Birkh\"auser, 1999, xv+510 p.
\bibitem{Hel}{\it A. Heller}, Stable homotopy categories, Bull.
         Amer. Math. Soc. 74 (1968), 28-63.
\bibitem{Hir}{\it Ph. S. Hirschhorn}, Model categories and their
         localizations, Mathematical Surveys and Monographs 99, 2003, xv+457 p.
\bibitem{H}{\it M. Hovey}, Model categories, Mathematical Surveys and
         Monographs 63, 1999, xii+209 p.
\bibitem{KV}{\it M. Karoubi, O. Villamayor}, Foncteurs $K^n$ en
         alg\`ebre et en topologie, C.~R. Acad. Sci. Paris 269 (1969), 416-419.
\bibitem{MV} {\it F. Morel, V. Voevodsky}, $\bb A^1$-homotopy theory of schemes,
         Publ. Math. l'I.H.E.S 90 (1999), 45-143.
\bibitem{Q}{\it D.~Quillen}, Homotopical algebra, Lecture
         Notes in Mathematics, No.~43, Springer-Verlag, 1967.
\bibitem{Q1}{\it D.~Quillen}, Higher algebraic K-theory.~I, In Algebraic K-theory~I, Lecture
         Notes in Mathematics, No.~341, Springer-Verlag, 1973, pp.~85-147.
\bibitem{S}{\it R. G. Swan}, Excision in algebraic K-theory. J. Pure Appl.
         Algebra 1(3) (1971), 221--252.
\bibitem{T}{\it R. W. Thomason, T. Trobaugh}, Higher algebraic $K$-theory of schemes and of
         derived categories, The Grothendieck Festschrift~III, Collect.
         Artic. in Honor of the 60th Birthday of A. Grothendieck, Progress in
         Mathematics 88, Birkh\"auser, 1990, pp.~247-435.
\bibitem{Ver}{\it J.-L. Verdier}, Des cat\'egories d\'eriv\'ees des cat\'egories ab\'eliennes. Ast\'erisque
         239 (1996), ix+253 pp.
\bibitem{V}{\it V. Voevodsky}, Homotopy theory of simplicial sheaves in completely decomposable
         topologies, K-theory Preprint Archives 443 (2000).
\bibitem{Wal}{\it F. Waldhausen}, Algebraic K-theory of spaces, In Algebraic and
         geometric topology, Proc. Conf., New Brunswick/USA 1983, Lecture
         Notes in Mathematics, No.~1126, Springer-Verlag, 1985, pp.~318-419.
\bibitem{W2}{\it C. Weibel}, KV-theory of categories, Trans. Amer.
         Math. Soc. 267(2) (1981), 621-635.
\bibitem{W1}{\it C. Weibel}, Homotopy algebraic K-theory, Contemp. Math. 83 (1989), 461-488.
\bibitem{W}{\it C. Weibel}, An introduction to algebraic K-theory,
         an electronic book in progress (available at {\tt
         math.rutgers.edu/$\sim$weibel}).
\end{thebibliography}
\end{document}